\documentclass[a4paper,10pt]{amsart}
\usepackage{amssymb, amsmath,latexsym,amsfonts,amsbsy, amsthm,mathtools,graphicx,CJKutf8,CJKnumb,CJKulem,color,geometry}
\usepackage{verbatim}
\def\l{\left}
\def\r{\right}
\def\f{\frac}

\def\az{\alpha}
\def\lz{\lambda}
\def\Lz{\Lambda}
\def\az{\alpha}
\def\rz{\rho}
\def\ez{\epsilon}
\def\bz{\beta}
\def\dz{\delta}
\def\gz{\gamma}
\def\tz{\theta}
\def\sz{\sigma}

\def\beq{\begin{equation}}
\def\beq*{\begin{equation*}}

\author{Ovidiu Savin and Qian Zhang}
\title{\bf{Boundary regularity for Monge-Amp\`ere equations with unbounded right hand side}}
\date{}

\theoremstyle{plain}
\theoremstyle{plain}\newtheorem{thm}{Theorem}[section]
\theoremstyle{plain}\newtheorem{prop}{Proposition}[section]
\theoremstyle{plain}
\theoremstyle{plain}\newtheorem{lem}{Lemma}[section]
\theoremstyle{plain}

\numberwithin{equation}{section}

\begin{document}
\maketitle
\noindent {\bf{Abstract}.}\quad
We consider Monge-Amp\`ere equations with right hand side $f$ that degenerate to $\infty$ near the boundary of a convex domain $\Omega$, which are of the type
$$\mathrm{det}\;D^2 u=f\quad\mathrm{in}\;\Omega,\quad\quad f\sim d^{-\az}_{\partial\Omega}\quad\mathrm{near}\;\partial\Omega,$$
where $d_{\partial\Omega}$ represents the distance to $\partial \Omega$ and $-\az$ is a negative power with $\az\in(0,2)$. We study the boundary regularity of the solutions and establish a localization theorem for boundary sections.
\bigskip

\section{Introduction}\label{s1}

In this paper we consider degenerate Monge-Amp\`ere equations of the type
\begin{equation}\label{s1: eq 1}
\mathrm{det}\;D^2 u=f\quad\mathrm{in}\;\Omega,\quad\quad f\sim d^{-\az}_{\partial\Omega}\quad\mathrm{near}\;\partial\Omega,
\end{equation}
where $d_{\partial\Omega}$ represents the distance to the boundary of the domain $\Omega$ and $-\az$ is a negative power with $\az\in(0,2)$. 

Boundary estimates for the Monge-Amp\`ere equation in the nondegenerate case $f\in C^2(\overline{\Omega})$, $f>0$, were obtained by Ivo$\check{c}$kina \cite{I}, Krylov \cite{K}, Caffarelli-Nirenberg-Spruck \cite{CNS} (see also \cite{C, TW}).

In \cite{S2}, a localization theorem at boundary points was proved when the right hand side $f$ is only bounded away from $0$ and $\infty$. It states that under natural local assumptions on the domain and boundary data, the sections $S_h(x_0)$ with $x_0\in\partial\Omega$ are ``equivalent" to half-ellipsoids centered at $x_0$. This extends up to the boundary a result that is valid for sections compactly included in $\Omega$, which is a consequence of John's lemma from convex geometry. These localization theorems are the key ingredients in establishing optimal $C^{2,\az}$ and $W^{2,p}$ estimates for solutions under further regularity properties of the right-hand side $f$ and boundary data (see \cite{C1, S2, S4}). 

In \cite{S3}, the first author studied degenerate Monge-Amp\`ere equations of the type
\begin{equation}\label{s1: eq 2}
\mathrm{det}\;D^2 u=f\quad\mathrm{in}\;\Omega,\quad\quad f\sim d^{\az}_{\partial\Omega}\quad\mathrm{near}\;\partial\Omega,
\end{equation}
where $\az>0$ is a positive power. A localization theorem and pointwise $C^2$ estimate were established in \cite{S3} and they were later used in \cite{LS} to prove the global smoothness for the eigenfunctions of the Monge-Amp\`ere operator $(\mathrm{det}\;D^2 u)^{1/n}$.

In this paper, we consider the case of the Monge-Amp\`ere equation with right hand side which degenerates to $\infty$ near the boundary of $\Omega$. This type of equations appear for example in the study of affine spheres in gemetry \cite{Ca,CY}, the
$p$-Minkowski problem \cite{Lu}, or in optimal transportation problems involving two densities with only one of them having compact support.

We study the case when $f$ is ``comparable" with a negative power $d_{\partial \Omega}^{-\alpha}$ of the distance function to $\partial \Omega$. It can be checked from a simple 1D example that the Dirichlet problem for equation \eqref{s1: eq 1} is well posed only for $\alpha \in (0,2)$. Moreover, when $\alpha \in (0,1)$ solutions are expected to have bounded gradients, and when $\alpha \in [1,2)$ the gradient should tend to $\infty$ as we approach the boundary. We study the geometry of boundary sections of solutions to \eqref{s1: eq 1} and prove two localization theorems Theorems \ref{s2: thm 2.1} and \ref{s1: thm 4} depending whether $\alpha$ is smaller or larger than 1. 

We first give the localization theorem for the case $\az\in(0,1)$. It states that under appropriate assumptions on the domain and boundary data, the sections 
$$S_{h}(x_0):=\{x\in\bar{\Omega}|\; u(x)<u(x_0)+\nabla u(x_0)\cdot(x-x_0)+h\}$$
with $x_0\in\partial\Omega$ have the shape of half-ellipsoids centered at $x_0$. 

\begin{thm}\label{s2: thm 2.1}
Assume $\Omega\subset\mathbb{R}^n$ is a bounded convex set, $\partial\Omega\in C^2$. Let $u:\overline{\Omega}\to\mathbb{R}$ be continuous, convex, satisfying
\begin{equation}\label{s1: eq 3.1}
\mathrm{det}\;D^2 u=f,\quad\quad \lz_0 d_{\partial\Omega}^{-\az}\le f\le\Lz_0d_{\partial\Omega}^{-\az}\quad\mathrm{in}\;\Omega
\end{equation}
for some $\az\in(0,1)$, and on $\partial\Omega$, $u$ separates quadratically from its tangent plane, namely
\begin{equation}\label{s1: eq 3.4}
\mu|x-x_0|^2\le u(x)-u(x_0)-\nabla u(x_0)\cdot(x-x_0)\le\mu^{-1}|x-x_0|^2,\quad\forall x, x_0\in\partial\Omega,
\end{equation}
for some $\mu>0$. Then there is a constant $c>0$ depending only on $n,\lz_0,\Lz_0,\az,\mu,\mathrm{diam}(\Omega)$ and $\|\partial\Omega\|_{C^2}$ such that for each $x_0\in\partial\Omega$ and $h\le c$ we have
\begin{equation*}
\mathcal{E}_{ch}(x_0)\cap\overline{\Omega}\subset S_h(x_0)\subset\mathcal{E}_{c^{-1}h}(x_0),
\end{equation*}
where
$$\mathcal{E}_h(x_0):=\{|(x-x_0)_\tau|^2+|(x-x_0)\cdot\nu_{x_0}|^{2-\az}<h\},\quad\forall h>0,$$
$\nu_{x_0}$ denotes the unit inner normal to $\partial\Omega$ at $x_0$ and
$$(x-x_0)_\tau:=(x-x_0)-[(x-x_0)\cdot\nu_{x_0}]\nu_{x_0}$$
is the projection of $x-x_0$ onto the tangent plane of $\partial\Omega$ at $x_0$.
\end{thm}

Theorem \ref{s2: thm 2.1} states that a boundary section $S_h$ is equivalent to an ellipsoid of axes $h^{\f{1}{2}}$ in the tangential direction to $\partial\Omega$ and $h^{\f{1}{2-\az}}$ in the normal. As a corollary, it can be proved that the maximal interior sections have the same geometry as boundary sections. Namely, for any $y_0\in\Omega$, let $S_{\bar{h}}(y_0)$ denote the maximal interior section centered at $y_0$ which becomes tangent to $\partial\Omega$ at some point $x_0$. Then $S_{\bar{h}}(y_0)$ is equivalent to an ellipsoid of axes $\bar{h}^{\f{1}{2}}$ in the tangential direction to $\partial\Omega$ at $x_0$ and $\bar{h}^{\f{1}{2-\az}}$ in the normal $\nu_{x_0}$.

We remark that if $u|_{\partial\Omega}=\varphi$ and $\partial\Omega\in C^3,\varphi\in C^3(\partial\Omega)$, and $\Omega$ is uniformly convex, then the quadratic separation condition \eqref{s1: eq 3.4} is satisfied. The proof is given in \cite[Proposition 3.2]{S2}, where only the lower bound of $\mathrm{det}\;D^2 u$ is used. Since in our degenerate case, $\mathrm{det}\;D^2 u$ is also bounded below by a constant, the estimate still applies.

Theorem \ref{s2: thm 2.1} implies global $W^{2,p}$ estimates of solutions if we assume further that $f=g \, d_{\partial\Omega}^{-\az}$ for some function $g\in C(\overline{\Omega})$ which is strictly positive. In a subsequent work we will show that $u\in W^{2,p}(\Omega)$ for any $p<\f{1}{\az}$.

For the case $\az\in(0,1)$, we establish the following Liouville type theorem for global solutions to \eqref{s1: eq 1}.

\begin{thm}\label{s1: thm 2}
Let $u\in C(\overline{\mathbb{R}^n_+})$ be a convex function that satisfies 
\begin{equation}\label{s5: eq 7.3}
c_0(|x'|^2+x_n^{2-\az})\le u(x)\le c_0^{-1}(|x'|^2+x_n^{2-\az})
\end{equation}
for some $c_0>0$ and
\begin{equation}\label{s5: eq 7.4}
\mathrm{det}\;D^2 u=x_n^{-\az},\quad\quad u(x',0)=\f{1}{2}|x'|^2.
\end{equation}
Then
\begin{equation*}
u(x)=\f{1}{2}|x'|^2+\f{x_n^{2-\az}}{(2-\az)(1-\az)}.
\end{equation*}
\end{thm}

Theorem \ref{s2: thm 2.1} and the Liouville theorem imply a pointwise $C^2$ tangential estimate at the boundary.

\begin{thm}\label{s1: thm 3}
Assume that $\Omega\subset\{x_n>0\}$ is a bounded convex set, $0\in\partial\Omega$, $\partial\Omega\in C^2$ near the origin, and the principal curvatures of $\partial\Omega$ at 0 are strictly positive. Let $u\in C(\overline{\Omega})$ be a convex solution to the equation
\begin{equation*}
\mathrm{det}\;D^2 u=f(x)d_{\partial\Omega}^{-\az}\quad\mathrm{in}\;\Omega,\quad\quad u=\varphi\quad\mathrm{on}\;\partial\Omega.
\end{equation*}
for some $\az\in(0,1)$, where $f$ is a nonnegative function that is continuous at the origin and $f(0)>0$, the boundary data $\varphi$ is $C^2$ at $0$, and it separates quadratically away from $0$. Assume further that
$$u(0)=0,\quad\nabla u(0)=0.$$
Then there exists a constant $a>0$ such that
\begin{equation*}
u(x)=Q(x')+a x_n^{2-\az}+o(|x'|^2+x_n^{2-\az}),
\end{equation*}
where $Q$ represents the quadratic part of the boundary data $\varphi$ at the origin.
\end{thm}

Next we give the localization theorem when $\az\in(1,2)$. In this case we consider the maximal sections included in $\Omega$ which become tangent to $\partial\Omega$ at boundary points.

\begin{thm}\label{s1: thm 4}
Assume $\Omega\subset\mathbb{R}^n$ is uniformly convex, $\partial\Omega\in C^2$. Assume further that $0\in\partial\Omega$ and the $x_n$ coordinate axis lies in the direction $\nu_{0}$ ($\nu_{0}$ is the unit inner normal to $\partial\Omega$ at $0$).

 Let $u:\overline{\Omega}\to\mathbb{R}$ be continuous, convex, satisfying
$$
\mathrm{det}\;D^2 u=f,\quad\quad\lz_0d_{\partial\Omega}^{-\az}\le f\le\Lz_0d_{\partial\Omega}^{-\az}\quad\mathrm{in}\;\Omega,$$
for some $\az\in[1,2)$, and assume $u|_{\partial\Omega}=\varphi \in C^2$. Suppose that $S_{\bar{h}}(y_0)$ is the maximal section included in $\Omega$ which becomes tangent to $\partial\Omega$ at $0$. Then 
$$\nabla_{x'}u(y_0)=\nabla_{x'}\varphi(0),\quad\quad M=-u_n(y_0)\ge -C,$$
and the following hold:
\begin{enumerate}
\item[i)] If $\az\in(1,2)$, denote $\bz:=\f{n+\az-1}{n}$, then we have
$$c\bar{h}^{\f{1-\bz}{2-\bz}}\le\max\{M,1\}\le C\bar{h}^{\f{1-\bz}{2-\bz}},\quad\quad c\bar{h}^{\f{1}{2-\bz}}\le d_{\partial\Omega}(y_0)\le C\bar{h}^{\f{1}{2-\bz}},$$
\begin{eqnarray*}
\{|x'|^{2}+|x_n|\le c\bar{h}^{\f{1}{2-\bz}}\}\subset S_{\bar{h}}(y_0)-y_0\subset\{|x'|^{2}+|x_n|\le C\bar{h}^{\f{1}{2-\bz}}\}.
\end{eqnarray*}
\item[ii)] If $\az=1$, denote $\bar{h}_*:=\min\{\bar{h},1\}$, then we have
$$-c\log(C\bar{h})\le|M|^n\le-C\log(c\bar{h}),\quad\quad c\bar{h}_*^{C}\le d_{\partial\Omega}(y_0)\le C\bar{h}_*^c,$$
\begin{eqnarray*}
B_{c\bar{h}_*^{C}}\subset&S_{\bar{h}}(y_0)-y_0\subset B_{C\bar{h}_*^c}.
\end{eqnarray*}
\end{enumerate}
Here the constants $c, C$ depend only on $n,\lz_0,\Lz_0,\az,\mathrm{diam}(\Omega)$, and $\varphi,\partial\Omega$ up to their second derivatives.
\end{thm}

In the case $\az\in(1,2)$, Theorem \ref{s1: thm 4} states that for any $y_0\in\Omega$, the maximal interior section $S_{\bar{h}}(y_0)$ which becomes tangent to $\partial\Omega$ at some point $x_0$ is equivalent to an ellipsoid of axes $\bar{h}^{\f{1}{2(2-\bz)}}$ in the tangential direction to $\partial\Omega$ at $x_0$ and $\bar{h}^{\f{1}{2-\bz}}$ in the normal $\nu_{x_0}$.
For the border line case $\alpha=1$, it cannot be concluded from ii) that $S_{\bar{h}}(y_0)$ is equivalent to an ellipsoid whose shape depends only on $\bar h, y_0$ and $\Omega$. Probably more precise information is needed on the ratio between $f$ and $d_{\partial \Omega}^{-1}$ in order to reach a similar conclusion as in the case $\alpha \in (1,2)$.  

The proofs of Theorems \ref{s2: thm 2.1} and \ref{s1: thm 4} are quite different. Theorem \ref{s1: thm 4} follows directly from comparison with explicit barriers. Theorem \ref{s2: thm 2.1} is much more involved and most of the paper will be devoted towards its proof. We will follow similar ideas as in the nondegenerate case treated in \cite{S2}.

The paper is organized as follows. In Section 2 we introduce some notation, then we reduce Theorem \ref{s2: thm 2.1} to its local version Theorem \ref{s4: thm 2.1}. This is further reduced to Theorem \ref{s4: thm 2.1f}, where the distance function is replaced by $x_n$. We also give a more precise quantitative version of Theorem \ref{s1: thm 3} (see Theorem \ref{s6: thm 2.4}). Sections 3-4 are devoted to the proof of Theorem \ref{s4: thm 2.1f}. In Section 5, the proof of Theorem \ref{s1: thm 2} is given. In Section 6, we give the proof of Theorem \ref{s6: thm 2.4} and then finish the proof of Theorem \ref{s1: thm 3}. In the last section, we give the proof of Theorem \ref{s1: thm 4}.

\section{Statement of main results}\label{s2}
We introduce some notation. We denote points in $\mathbb{R}^n$ as 
$$x=(x_1,\dots,x_n)=(x',x_n),\quad x'\in\mathbb{R}^{n-1}.$$ 
Let $u$ be a convex function defined on a convex set $\overline{\Omega}$, we denote by $S_h(x_0)$ the section centered at $x_0$ and at height $h>0$,
$$S_{h}(x_0):=\{x\in\overline{\Omega}|\; u(x)<u(x_0)+\nabla u(x_0)\cdot(x-x_0)+h\}.$$
When $x_0\in\partial\Omega$, the term $\nabla u(x_0)$ is understood in the sense that
$$x_{n+1}=u(x_0)+\nabla u(x_0)\cdot(x-x_0)$$
is a supporting hyperplane for the graph of $u$ at $x_0$ but for any $\ez>0$,
$$x_{n+1}=u(x_0)+(\nabla u(x_0)+\ez\nu_{x_0})\cdot(x-x_0)$$
is not a supporting hyperplane, where $\nu_{x_0}$ denotes the unit inner normal to $\partial\Omega$ at $x_0$. We denote for simplicity $S_{h}=S_{h}(0)$, and sometimes when we specify the dependence on the function $u$ we use the notation $S_{h}(u)=S_{h}$.

For a set $E\subset\mathbb{R}^n$, we always denote $\pi(E)$ the projection of $E$ into $\mathbb{R}^{n-1}$, i.e., 
$$\pi(E):=\{x'\in\mathbb{R}^{n-1}: \exists\;t\in\mathbb{R}\;s.t.\;(x', t)\in E\}.$$
In the case $\az\in(0,1)$, for any $h>0$ we often use the particular sets
$$\mathcal{E}_h:=\{|x'|^2+x_n^{2-\az}<h\},\quad\quad\mathcal{E}_h^+:=\mathcal{E}_h\cap\{x_n>0\},$$
and the diagonal matrix
$$F_h:=\mathrm{diag}(h^{\f{1}{2}},h^{\f{1}{2}},\dots,h^{\f{1}{2}},h^{\f{1}{2-\az}})$$
in our estimates.

Next we give a local version of Theorem \ref{s2: thm 2.1}. Our assumptions are the following.

Let $\Omega\subset\mathbb{R}^n$ be a open convex set. Assume that for some fixed small $\rz>0$,
\begin{equation}\label{s4: eq 3.1}
B_\rz(\rz e_n)\subset\Omega\subset\{x_n>0\}\cap B_{\f{1}{\rz}},
\end{equation}
and 
\begin{equation}\label{s4: eq 3.1'}
\Omega\;\mathrm{contains\;an\;interior\;ball\;of\;radius}\;\rz\;\mathrm{tangent\;to}\;\partial\Omega\;\mathrm{at\;each\;point\;on}\;\partial\Omega\cap\{x_n\le\rz\}.
\end{equation}
The part $\partial\Omega\cap\{x_n\le\rz\}$ is then given by $x_n=g(x')$ for some convex function $g$, where
\begin{equation}\label{s4: eq 3.1''}
g\in C^2\l(\pi(\partial\Omega\cap\{x_n<\rz\})\r), \quad g(0)=0, \quad\nabla g(0)=0.
\end{equation}

Let $u:\overline{\Omega}\to\mathbb{R}$ be a convex solution to
\begin{equation}\label{s4: eq 3.2}
\mathrm{det}\;D^2 u=f,\quad 0<\lz(x_n-g)^{-\az}\le f\le\Lz(x_n-g)^{-\az}\quad\mathrm{in}\;\Omega\cap\{x_n<\rz/2\}
\end{equation}
for some $\az\in(0,1)$. Moreover,
\begin{equation}\label{s4: eq 3.3}
x_{n+1}=0\;\mathrm{is\;the\;tangent\;plane\;to}\;u\;\mathrm{at}\;0,
\end{equation}
that is,

$u\ge 0,\;u(0)=0,\;\nabla u(0)=0$ in the sense that $x_{n+1}=\ez x_n$ is not a supporting plane for the graph of $u$ at $0$ for any $\ez>0$.

We also assume that $u$ separates quadratically on $\partial\Omega$ (in a neighborhood of $\{x_n=0\}$) from the tangent plane at $0$, i.e.,
\begin{equation}\label{s4: eq 3.4}
\mu |x|^2\le u(x)\le\mu^{-1}|x|^2\quad\mathrm{on}\;\partial\Omega\cap\{x_n\le\rz\}.
\end{equation}

\begin{thm}\label{s4: thm 2.1}
Assume $\Omega$ and $u$ satisfy \eqref{s4: eq 3.1}-\eqref{s4: eq 3.4}. Then there is a constant $c>0$ depending only on $n,\lz,\Lz,\az,\mu$ and $\rz$ such that for each $h\le c$ we have
\begin{equation*}
\mathcal{E}_{ch}\cap\overline{\Omega}\subset S_h\subset\mathcal{E}_{c^{-1}h}.
\end{equation*}
\end{thm}

Assume $\Omega$ and $u$ satisfy the hypotheses in Theorem \ref{s2: thm 2.1}. Fix a point $x_0\in\partial\Omega$, by a translation and a rotation of coordinates we can assume that $x_0=0$, and the $x_n$ coordinate axis lies in the direction $\nu_{x_0}$. Since $\partial\Omega\in C^2$, there exists $\rz>0$ such that \eqref{s4: eq 3.1}-\eqref{s4: eq 3.1''} hold, and after subtracting a linear function we have \eqref{s4: eq 3.3} and \eqref{s4: eq 3.4}. By \eqref{s4: eq 3.1}-\eqref{s4: eq 3.1''}, it is easy to see that
\begin{equation}\label{s4: eq k6}
\|D^2 g\|_{C\l(\pi(\partial\Omega\cap\{x_n\le\rz/2\})\r)}\le C(n,\rz)
\end{equation} 
and therefore
$$d_{\partial\Omega}(x)\le x_n-g(x')\le C'(n,\rz)d_{\partial\Omega}(x),$$
where $C(n,\rz)$ and $C'(n,\rz)$ are constants depending only on $n$ and $\rz$. It follows that $u$ satisfies \eqref{s4: eq 3.2} with $\lz:=\lz_0, \Lz:=C'(n,\rz)\Lz_0$. Therefore we reduce the proof of Theorem \ref{s2: thm 2.1} to that of Theorem \ref{s4: thm 2.1} above.

Let $\Omega$ and $u$ satisfy the hypotheses in Theorem \ref{s4: thm 2.1}. By constructing some lower barrier for $u$, we will prove in Section \ref{s3} that in some domain $\Omega_0\subset\Omega$ we have $x_n-g\sim x_n$, and $u$ still satisfies the quadratic separation \eqref{s4: eq 3.4} on $\partial\Omega_0$ in a neighborhood of $\{x_n=0\}$. Therefore we reduce the proof of Theorem \ref{s4: thm 2.1} to that of Theorem \ref{s4: thm 2.1f} below.

We assume \eqref{s4: eq 3.1}, \eqref{s4: eq 3.3}, \eqref{s4: eq 3.4} hold while replacing the equation \eqref{s4: eq 3.2} by 
\begin{equation}\label{s4: eq 3.2f}
\mathrm{det}\;D^2 u=f,\quad 0<\lz x_n^{-\az}\le f\le\Lz x_n^{-\az}\quad\mathrm{in}\;\Omega\cap\{x_n<\rz\}.
\end{equation}
Note that we do not assume \eqref{s4: eq 3.1'} and \eqref{s4: eq 3.1''} hold here.

\begin{thm}\label{s4: thm 2.1f}
Assume $\Omega$ and $u$ satisfy \eqref{s4: eq 3.1}, \eqref{s4: eq 3.3}, \eqref{s4: eq 3.4} and \eqref{s4: eq 3.2f}. Then there is a constant $c>0$ depending only on $n,\lz,\Lz,\az,\mu$ and $\rz$ such that for each $h\le c$ we have
\begin{equation*}
\mathcal{E}_{ch}\cap\overline{\Omega}\subset S_h\subset\mathcal{E}_{c^{-1}h}.
\end{equation*}
\end{thm}

We prove Theorem \ref{s4: thm 2.1f} using the compactness methods in \cite{S2}. We first obtain some preliminary estimates about $u$. Next we consider the rescaling $v$ of $u$. Then we reduce the proof of the theorem to that of a statement about $v$. We reduce this to the proof of a statement (Proposition \ref{s4: prop 5.2'}) about the limiting function (still denoted by $u$) of such $v$. Different from the case that $\az=0$ (in this case the estimate of the volume of $S_t(v)$ is $|S_t(v)|^2\sim t^n$), the estimate of the volume of $S_t(v)$ becomes
\begin{equation}\label{s2: eq 1}
(x^*_t(v)\cdot e_n)^{-\az}|S_t(v)|^2\sim t^n,
\end{equation}
where $x^*_t(v)$ is the center of mass of $S_t(v)$. The limiting function $u$ also satisfies this estimate. To prove Proposition \ref{s4: prop 5.2'}, we construct some lower barrier for the limiting function $u$ and use \eqref{s2: eq 1}. Since we do not have the estimate of $|S_h(u)|$, we also use the convexity of the original solution to estimate the quantity $x^*_t(v)\cdot e_n$ from below. The estimate \eqref{s2: eq 1} brings another difficulty when we prove Proposition \ref{s4: prop 5.2'}. We use John's lemma and find an ellipsoid $E_h$ equivalent to the section $S_h(u)$ of the limiting solution $u$. In the case $\az=0$, we use the estimate $|E_h|^2\sim h^n$ to estimate the shape of $S_h(u)$, but in our degenerate case, we do not have the estimate of the volume of $E_h$. For this, we use the estimate $(x^*_h(u)\cdot e_n)^{-\az}|E_h|^2\sim h^n$ to obtain an estimate of the shape of $S_h(u)$ in terms of the quantity $x^*_h(u)\cdot e_n$. Using this estimate, we rescale $u$ and reduce Proposition \ref{s4: prop 5.2'} to the lower-dimensional case. Again, since we do not have the estimate of $|E_h|$, we perform a different rescaling (which corresponds to our estimate \eqref{s2: eq 1}) from the $\az=0$ case.

At the end of this section we give a more precise quantitative version of Theorem \ref{s1: thm 3} as follows.

\begin{thm}\label{s6: thm 2.4}
For any $\eta>0$ there exists $\ez_0>0$ depending only on $\eta, n, \az$ such that if \eqref{s4: eq 3.1}-\eqref{s4: eq 3.3} hold with $\lz=1-\ez_0, \Lz=1+\ez_0$ and 
\begin{equation}\label{s6: thm 2.44}
\l(\f{1}{2}-\ez_0\r)|x'|^2\le u(x)\le\l(\f{1}{2}+\ez_0\r)|x'|^2\quad\mathrm{on}\;\partial\Omega\cap\{x_n\le\rz\},
\end{equation}
then for all $h\le c$, we have
\begin{equation*}
(1-\eta)S_h(U_0)\cap\overline{\Omega}\subset S_h(u)\subset(1+\eta)S_h(U_0),
\end{equation*}
where 
$$U_0(x):=\f{1}{2}|x'|^2+\f{x_n^{2-\az}}{(2-\az)(1-\az)},\quad\quad S_h(U_0):=\{x\in\mathbb{R}^n: U_0(x)<h\},$$
and the constant $c>0$ depends only on $\eta, n,\az,\rz$.
\end{thm}

\section{Proof of Theorem \ref{s4: thm 2.1f} (I)}\label{s3}

As mentioned in Section \ref{s2}, we first show that we can reduce the proof of Theorem \ref{s4: thm 2.1} to that of Theorem \ref{s4: thm 2.1f}.

\begin{prop}\label{s4: prop1}
Theorem \ref{s4: thm 2.1f} implies Theorem \ref{s4: thm 2.1}.
\begin{proof}
In this proof we always denote by $c, C, c_i, C_i(i=0,1,2,\dots)$ constants depending only on $n,\lz,\Lz,\mu,\az$ and $\rz$. For simplicity of notation, their values may change from line to line whenever there is no possibility of confusion.

Let
$$v_0:=\mu|x'|^2+\f{\Lz}{(2-\az)(1-\az)\,\mu^{n-1}}(x_n-g)^{2-\az}.$$
Then by straightforward computation and using \eqref{s4: eq k6}, we obtain that
\begin{eqnarray}\label{v_0}
\mathrm{det}\;D^2 v_0&=&\f{\Lz}{\mu^{n-1}}(x_n-g)^{-\az}\mathrm{det}\l(2\mu I_{n-1}-\f{\Lz}{(1-\az)\,\mu^{n-1}}(x_n-g)^{1-\az}D^2 g\r)\nonumber\\
&\ge&\Lz (x_n-g)^{-\az}
\quad\mathrm{in}\;\Omega\cap\{x_n<c_*\},
\end{eqnarray}
where $c_*\le\rz/2$ is small depending only on $n,\Lz,\mu,\az$ and $\rz$.

Denote $D:=\pi\l(\Omega\cap\{x_n=c_*\}\r)$. For $x'\in D$, define
$$g^*(x'):=\sup\l\{l(x'): l\le g\;\mathrm{in}\;D,\;l\;\mathrm{is\;affine,\;and}\;|\nabla l|\le\f{c_*\rz}{2}\r\}.$$
Then $g^*$ is convex in $D$ since it is the supremum of a family of convex functions.

We claim that for any $x\in\Omega\cap\{x_n=c_*\}$, we have
\begin{eqnarray}\label{s4: eq g^*1}
x_n-g^*(x')\ge\f{c_*}{2}.
\end{eqnarray}
Indeed, if $l$ is affine, $l\le g$ in $D$ and $|\nabla l|\le\f{c_*\rz}{2}$, then 
$$0=g(0)\ge l(0)=l(x')-\nabla l\cdot x',$$
it follows that
$$l(x')\le\nabla l\cdot x'\le\f{c_*\rz}{2}\cdot\f{1}{\rz}=\f{c_*}{2},$$
where we use the fact that $\Omega\subset B^+_{1/\rz}$. Thus the claim follows. 

We also claim that
\begin{equation}\label{s4: eq g^*1'}
\pi(\Omega\cap\{x_n\le c_0\rz\})\subset D\cap\l\{|\nabla g|\le\f{c_*\rz}{2}\r\}\subset\{g^*=g\}
\end{equation}
for some small constant $c_0$. 

Indeed, the second inclusion in \eqref{s4: eq g^*1'} follows easily from the convexity of $g$ and the definition of $g^*$. Therefore we only need to prove the first inclusion. Let $c_0>0$ be a small constant to be chosen. For any $x_0\in\partial\Omega\cap\{x_n\le c_0\rz\}$, we have $B_\rz(y_0)\subset\Omega\subset\{x_n\ge 0\}$ by \eqref{s4: eq 3.1'}, where $y_0:=x_0+\rz\nu_{x_0}$. Let 
$$t=\inf_{x\in B_\rz(y_0)}x_n,$$
then $(y'_0, t)\in\partial B_\rz(y_0)$ and 
$$\rz\nu_{x_0}\cdot e_n=(y_0\cdot e_n-t)-(x_0\cdot e_n-t)=\rz-(x_0\cdot e_n-t)\ge (1-c_0)\rz,$$
which gives
$$\f{1}{\sqrt{1+|\nabla g(x'_0)|^2}}=\nu_{x_0}\cdot e_n\ge 1-c_0.$$
Hence,
\begin{equation}\label{s4: eq k5}
|\nabla g(x'_0)|\le\sqrt{\l(\f{1}{1-c_0}\r)^2-1}\le\f{c_*\rz}{2}
\end{equation}
if $c_0$ is small. The desired conclusion \eqref{s4: eq g^*1'} follows.

Let
$$v^*:=\mu|x'|^2+\f{\Lz}{(2-\az)(1-\az)\,\mu^{n-1}}(x_n-g)^{2-\az}-C^*(x_n-g^*(x')).$$
Then $v^*$ is a lower barrier for $u$ in $\Omega\cap\{x_n\le c_*\}$ if $C^*$ is large depending only on $n,\Lz,\mu,\az$ and $\rz$. 

Indeed, since $g^*$ is convex, we find from \eqref{v_0} that $v^*$ is a subsolution of the equation
$$\mathrm{det}\;D^2 w=\Lz (x_n-g)^{-\az}.$$
On $\partial\Omega\cap\{x_n\le c_*\}$, we have $x_n-g^*=g-g^*\ge 0$, which implies
$$v^*\le\mu|x'|^2\le u.$$
On $\Omega\cap\{x_n=c_*\}$, we obtain from \eqref{s4: eq g^*1} that
$$v^*\le\f{\mu}{\rz^2}+\f{\Lz}{(2-\az)(1-\az)\,\mu^{n-1}} c_*^{2-\az}-C^*\f{c_*}{2}\le 0\le u$$
if $C^*$ is large.

Thus, 
$$v^*\le u\quad\mathrm{in}\;\Omega\cap\{x_n\le c_*\}.$$
This together with \eqref{s4: eq g^*1'} implies that
\begin{equation}\label{s4: eq 4.1f}
u\ge\mu|x'|^2-C^*(x_n-g(x'))\quad\mathrm{in}\;\Omega\cap\{x_n\le c_0\rz\}.
\end{equation}
Therefore, if $\dz$ is small, we have
\begin{equation}\label{s4: eq 3.4f0}
u\ge\f{\mu|x'|^2}{2}\quad\mathrm{in}\;\Omega\cap\{x_n\le c_0\rz\}\cap\{x_n\le g(x')+\dz|x'|^2\}.
\end{equation}

On the other hand, the convexity of $u$ and the quadratic separation of $u$ on $\partial\Omega\cap\{x_n\le\rz\}$ (see \eqref{s4: eq 3.4}) implies that
\begin{equation}\label{s4: eq 3.4f0'}
u\le C|x'|^2\quad\mathrm{in}\;\Omega\cap\{x_n\le c_0\rz\}\cap\{x_n\le g(x')+\dz|x'|^2\}.
\end{equation}
In particular, if we denote $\Omega_0:=\Omega\cap\{x_n<c_0\rz\}\cap\{x_n>g(x')+\dz|x'|^2\}$, then the above two estimates hold on $\partial\Omega_0\cap\{x_n\le c_0\rz\}$. 

We have
$$\lz x_n^{-\az}\le\mathrm{det}\;D^2 u\le Cx_n^{-\az}\quad\mathrm{in}\;\Omega_0.$$

We apply Theorem \ref{s4: thm 2.1f} to $u$ in $\Omega_0$ and obtain that
$$\mathcal{E}_{ch}\cap\overline{\Omega_0}\subset S_h\cap\overline{\Omega_0}\subset\mathcal{E}_{Ch},\quad\forall h\le c.$$

We claim that the last estimate also holds for $S_h$ (instead of $S_h\cap\overline{\Omega_0}$). Indeed, we have by \eqref{s4: eq 3.4f0'}
$$(\overline{\Omega}\setminus\Omega_0)\cap\mathcal{E}_{ch}\subset S_h$$
and therefore
$$\mathcal{E}_{ch}\cap\overline{\Omega}\subset S_h.$$
On the other hand, we obtain from \eqref{s4: eq 3.4f0} that
$$(S_h\cap\{x_n\le c_0\rz\})\setminus\Omega_0\subset\{|x'|\le Ch^{\f{1}{2}}\}.$$
Since 
$$(\overline{\Omega}\setminus\Omega_0)\cap\{x_n\le c_0\rz\}\subset\{x_n\le g(x')+\dz|x'|^2\},$$
we obtain
$$S_h\setminus\Omega_0\subset\{|x'|\le Ch^{\f{1}{2}}, x_n\le Ch\}\subset\mathcal{E}_{Ch}.$$
\end{proof}
\end{prop}

In the following we give the first part of the proof of Theorem \ref{s4: thm 2.1f}. In the remaining part of this section we denote by $c, C, c_i, C_i(i=0,1,2,\dots)$ positive constants depending on $n,\lz,\Lz,\mu$ and $\az$. The dependence of various constants also on $\rz$ will be denoted by $c(\rz), C(\rz), c_i(\rz), C_i(\rz)(i=0,1,2,\dots)$.
 
\begin{prop}\label{s4: prop 4.1f}
Assume that $\Omega$ and $u$ satisfy the hypotheses of Theorem \ref{s4: thm 2.1f}. Then, for each $h\le c(\rz)$ there exists a linear transformation (sliding along $x_n=0$)
$$A_h x=x-\nu x_n,\quad \nu_n=0,\quad |\nu|\le C(\rz)h^{-\f{n}{2(n+1-\az)}},$$
such that the rescaled function
$$\tilde{u}(A_h x)=u(x)$$
satisfies in
$$\tilde{S}_h:=A_h S_h=\{\tilde{u}<h\}$$
the following
\begin{enumerate}
\item[\rm(i)] the center of mass $\tilde{x}^*_h$ of $\tilde{S}_h$ lies on the $x_n$ axis, i.e. $\tilde{x}^*_h=d_h e_n$.
\item[\rm(ii)]
$$ch^n\le |S_h|^2 d_h^{-\az}\le Ch^n.$$
And after a rotation of the $x_1,\dots, x_{n-1}$ variables we have
$$\tilde{x}^*_h+cD_h B_1\subset \tilde{S}_h\subset CD_h B_1,$$
where $D_h:=\mathrm{diag}(d_1, d_2,\dots, d_{n-1}, d_n)$
is a diagonal matrix that satisfies
\begin{equation}\label{s4: eq d_nf}
\l(\prod_{1}^{n-1}d_i^2\r)d_n^{2-\az}=h^n
\end{equation}
and 
$$cd_h\le d_n\le Cd_h.$$
\item[\rm(iii)] Denote $\tilde{\Omega}_h:=A_h\Omega$ and $\tilde{G}_h:=\partial\tilde{S}_h\cap\{\tilde{u}<h\}$, then $\tilde{G}_h$ is a graph i.e.
$$\tilde{G}_h=(x',\tilde{g}_h(x'))\quad\quad\mathrm{with}\quad\quad\tilde{g}_h(x')\le\f{2}{\rz}|x'|^2$$
and the function $\tilde{u}$ satisfies on $\tilde{G}_h$
$$\f{\mu}{2}|x'|^2\le\tilde{u}(x)\le 2\mu^{-1}|x'|^2.$$
\end{enumerate}

\begin{proof}
Let
$$v:=\mu|x'|^2+\f{\Lz}{(2-\az)(1-\az)\,\mu^{n-1}}x_n^{2-\az}-C(\rz)x_n,$$
where $C(\rz)$ is large such that
$$\f{\Lz}{(2-\az)(1-\az)\,\mu^{n-1}}x_n^{2-\az}-\f{C(\rz)}{2}x_n\le 0\quad\mathrm{in}\;\Omega\cap\{x_n\le\rz\},$$
then it is straightforward to check that $v$ is a lower barrier for $u$ in $\Omega\cap\{x_n\le\rz\}$. It follows that
\begin{equation}\label{s4: eq 4.1f}
S_h\cap\{x_n\le\rz\}\subset\{v<h\}\subset\{x_n>c(\rz)(\mu|x'|^2-h)\}.
\end{equation}

Let $x^*_h$ be the center of mass of $S_h$ and $d_h:=x^*_h\cdot e_n$. We claim that
\begin{equation}\label{s4: eq 4.2f}
d_h\ge c_0(\rz) h^{\f{n}{n+1-\az}}
\end{equation}
for some small $c_0(\rz)>0$. 

Indeed, if 
$$d_h\ge c(n)\rz$$
with $c(n)$ depending only on $n$, then \eqref{s4: eq 4.2f} holds clearly. On the other hand, if 
$$d_h\le c(n)\rz,$$
then by John's lemma, for some constant $C(n)$ depending only on $n$ we have
$$S_h\subset\l\{x_n\le C(n)d_h\le\f{\rz}{2}\r\}$$
if $c(n)$ is small. If \eqref{s4: eq 4.2f} does not hold, then from the last estimate, \eqref{s4: eq 4.1f} and John's lemma that
$$S_h\subset\{x_n\le C(n)c_0(\rz)h^{\f{n}{n+1-\az}}\le h^{\f{n}{n+1-\az}}\}\cap\{|x'|\le C_1(\rz) h^{\f{n}{2(n+1-\az)}}\}.$$
Define
$$w=\ez x_n+\f{h}{2}\l(\f{|x'|}{C_1(\rz)h^{\f{n}{2(n+1-\az)}}}\r)^2+\f{\Lz [C_1(\rz)]^{2(n-1)}h}{(2-\az)(1-\az)}\l(\f{x_n}{h^{\f{n}{n+1-\az}}}\r)^{2-\az}.$$
Then we have in $S_h$,
$$w\le\ez+\f{h}{2}+\f{\Lz[C_1(\rz)]^{2(n-1)}h}{(2-\az)(1-\az)}[C(n)c_0(\rz)]^{2-\az}\le h$$
if $c_0(\rz)$ is small. On $S_h\cap\partial\Omega$,
$$w\le\f{\ez}{\rz}|x'|^2+\f{h^{\f{1-\az}{n+1-\az}}}{2C_1(\rz)^{2}}|x'|^2+\f{\Lz[C_1(\rz)]^{2(n-1)}h^{\f{1-\az}{n+1-\az}}}{(2-\az)(1-\az)}\cdot\f{|x'|^2}{\rz}\le\mu|x'|^2$$
if $h\le c(\rz)$. In conclusion, 
$$w\le u\quad\mathrm{in}\;S_h,$$
which contradicts that $\nabla u(0)=0$. Thus \eqref{s4: eq 4.2f} holds.

Now we prove that for all small $h$ we have
\begin{equation}\label{s4: eq *1f}
d_h\le C_0 h^{\f{1}{2-\az}}
\end{equation}
for some large constant $C_0$. 

Assume by contradiction that $d_h\ge C_0 h^{\f{1}{2-\az}}$. Then \eqref{s4: eq 4.1f} implies that
\begin{equation}\label{s4: eq *3f}
|(x^*_h)'|\le C(\rz) d_h^{\f{1}{2}}.
\end{equation}
From \eqref{s4: eq 3.1} and \eqref{s4: eq 3.4} we know that if $h\le c(\rz)$, then $S_h$ contains the set $\partial\Omega\cap\{x_n\le\rz\}\cap\{x: |x'|\le (ch)^{\f{1}{2}}\}$ for some small $c$ depending only on $\mu$. Therefore $S_h$ contains the convex set generated by $\partial\Omega\cap\{x_n\le\rz\}\cap\{x: |x'|\le (ch)^{\f{1}{2}}\}$ and the point $x^*_h$. Let $x_n=b$ be a hyperplane in $\mathbb{R}^n$, where $b\le\rz$ is chosen such that
$$ch+(b-\rz)^2=\rz^2.$$
For each $x_0\in\partial\Omega\cap\{x_n\le\rz\}\cap\{x: |x'|=(ch)^{\f{1}{2}}\}$, let $y_0$ be the intersection of the segment $\overline{x_0x^*_h}$ (which is the segment joining $x_0$ and $x^*_h$) and the hyperplane $x_n=b$. We can write 
$$y_0=(1-\tz)x_0+\tz x^*_h$$
for some $\tz=\tz(x_0)\in(0,1)$. Since 
$$(1-\tz)x_0\cdot e_n+\tz d_h=y_0\cdot e_n=b\le\f{ch}{\rz},$$
we obtain
\begin{equation*}
\tz\le\f{ch}{\rz d_h}.
\end{equation*}
Recall that $d_h\ge C_0 h^{\f{1}{2-\az}}$, then by \eqref{s4: eq *3f} we obtain that for all small $h$
\begin{eqnarray*}
|y'_0|=|(1-\tz)x'_0+\tz(x^*_h)'|&\ge&|x'_0|-\tz\l(|(x^*_h)'|+|x'_0|\r)\\
&\ge&(ch)^{\f{1}{2}}-\f{ch}{\rz d_h}\l(C(\rz)d_h^{\f{1}{2}}+(ch)^{\f{1}{2}}\r)\\
&\ge&\f{(ch)^{\f{1}{2}}}{2}.
\end{eqnarray*}
Since $S_h$ contains the convex set generated by all such $y_0$ and $x^*_h$, this means that $S_h$ contains a convex set of measure $c(n)\l(\f{(ch)^{\f{1}{2}}}{2}\r)^{n-1}d_h$, and therefore
\begin{equation}\label{s4: eq *f}
|S_h|\ge c(n)\l(\f{(ch)^{\f{1}{2}}}{2}\r)^{n-1}d_h.
\end{equation}
Let $v$ solves
$$\mathrm{det}\;D^2 v=\lz(C(n) d_h)^{-\az}\le\mathrm{det}\;D^2 u\quad\mathrm{in}\;S_h,\quad\quad v=h\quad\mathrm{on}\;\partial S_h.$$
Then
$$v\ge u\ge 0\quad\mathrm{in}\;S_h.$$
It follows 
$$h^n\ge|h-\min_{S_h}v|^n\ge c(n,\az)\lz d_h^{-\az}|S_h|^2.$$
Namely,
\begin{equation}\label{s4: eq 4.4f}
d_h^{-\az}|S_h|^2\le C(n,\lz,\az)h^n.
\end{equation}
It follows from \eqref{s4: eq *f} and \eqref{s4: eq 4.4f} that
$$d_h\le Ch^{\f{1}{2-\az}}.$$
We reach a contradiction if $C_0$ is sufficiently large, hence \eqref{s4: eq *1f} is proved.

Define 
$$A_h x=x-\nu x_n,\quad\nu=\f{(x^*_h)'}{d_h}$$
and
$$\tilde{u}(A_h x)=u(x).$$
Then the center of mass of $\tilde{S}_h=A_hS_h$ is
$$\tilde{x}^*_h=A_h x^*_h$$
and lies on the $x_n$-axis from the definition of $A_h$. We obtain from \eqref{s4: eq 4.1f} and \eqref{s4: eq 4.2f} that
\begin{equation}\label{s4: eq *10f}
|\nu|=\f{|(x^*_h)'|}{d_h}\le C(\rz)d_h^{-\f{ 1}{2}}\le C(\rz)h^{-\f{n}{2(n+1-\az)}}.
\end{equation}
Part $(i)$ of Proposition \ref{s4: prop 4.1f} follows.

Let $\tilde{\Omega}_h:=A_h\Omega$ and $\tilde{G}_h:=\partial\tilde{S}_h\cap\partial\tilde{\Omega}_h=\partial\tilde{S}_h\cap\{\tilde{u}<h\}$. 

On $\partial\Omega\cap\{x_n\le\rz\}\cap\{|x'|\le (\mu^{-1}h)^{\f{1}{2}}\}$, we have
$$|A_h x-x|=|\nu|x_n\le C(\rz)h^{-\f{n}{2(n+1-\az)}}|x'|^2\le C(\rz)h^{\f{1-\az}{2(n+1-\az)}}|x'|.$$
Note that $$\partial S_h\cap\partial\Omega\subset\{x_n\le\rz\}\cap\{|x'|\le (\mu^{-1}h)^{\f{1}{2}}\},$$
thus on $\tilde{G}_h=\partial\tilde{S}_h\cap\partial\tilde{\Omega}_h$,
$$x_n\le\f{1}{\rz}|(A_h^{-1} x)'|^2\le\f{2}{\rz}|x'|^2$$
and
$$\f{\mu}{2}|x'|^2\le\tilde{u}(x)=u(A_h^{-1}x)\le\mu^{-1}|(A_h^{-1}x)|^2\le 2\mu^{-1}|x'|^2.$$

It remains to prove $(ii)$. After a rotation of $x_1,\dots, x_{n-1}$ variables, we can assume that $\tilde{S}_h\cap\{x_n=d_h\}$ is equivalent to an ellipsoid of axes $d_1\le d_2\le\cdots\le d_{n-1}$ i.e.
$$\l\{\sum_1^{n-1}(\f{x_i}{d_i})^2\le 1\r\}\cap\{x_n=d_h\}\subset\tilde{S}_h\cap\{x_n=d_h\}\subset\l\{\sum_1^{n-1}(\f{x_i}{d_i})^2\le C(n)\r\}.$$
Thus,
$$\tilde{S}_h\subset\l\{\sum_1^{n-1}(\f{x_i}{d_i})^2\le C(n)\r\}\cap\{0\le x_n\le C(n)d_h\}.$$
Since $\tilde{u}\le 2\mu^{-1}|x'|^2$ on $\tilde{G}_h$, we see that the domain of definition of $\tilde{G}_h$ contains a ball in $\mathbb{R}^{n-1}$ of radius $(\mu h/2)^{\f{1}{2}}$. This implies that
\begin{equation}\label{s4: eq *4f}
d_i\ge c_1 h^{\f{1}{2}},\quad i=1,\dots, n-1.
\end{equation}

Now we prove that
\begin{equation}\label{s4: eq *5f}
d_h^{2-\az}\prod_1^{n-1}d_i^2\ge c_2 h^n.
\end{equation}
Indeed, if the last estimate does not hold, then we construct
$$w:=\ez x_n+\l[\sum_1^{n-1}(\f{x_i}{d_i})^2+(\f{x_n}{d_h})^{2-\az}\r]\cdot ch.$$
If $c_2$ is small, then we have
$$\mathrm{det}\;D^2 w\ge\f{c^n 2^{n-1}(2-\az)(1-\az)x_n^{-\az}}{c_2}>\Lz x_n^{-\az}.$$
On $\partial\tilde{S}_h\setminus\tilde{G}_h$,
$$w\le\ez+C(n,\az)ch\le h,$$
and on $\tilde{G}_h$, we use \eqref{s4: eq *4f} and \eqref{s4: eq 4.2f} to obtain
$$w\le \f{2\ez}{\rz}|x'|^2+ch\sum_1^{n-1}(\f{x_i}{d_i})^2+chC(n)^{1-\az}\f{2|x'|^2}{\rz d_h}\le\f{\mu}{2}|x'|^2$$
if $c$ is small. We conclude that $w\le\tilde{u}$ in $\tilde{S}_h$. This contradicts $\nabla\tilde{u}(0)=0$ and therefore \eqref{s4: eq *5f} holds.

Since $\tilde{S}_h$ contains the convex set generated by $\l\{\sum_1^{n-1}(\f{x_i}{d_i})^2\le 1\r\}\cap\{x_n=d_h\}$ and the point $0$, we have
$$|\tilde{S}_h|\ge c(n)\l(\prod_1^{n-1}d_i\r)\cdot d_h.$$
This together with \eqref{s4: eq *5f}, \eqref{s4: eq 4.4f} implies that
\begin{equation}\label{s4: eq *6f}
C h^n\ge d_h^{-\az}|\tilde{S}_h|^2\ge c(n)d_h^{2-\az}\prod_1^{n-1}d_i^2\ge ch^n.
\end{equation}
Define $d_n$ from $d_1, \dots, d_{n-1}$ by \eqref{s4: eq d_nf}, and \eqref{s4: eq *6f} gives
$$c d_h\le d_n\le C d_h.$$
This proves $(ii)$.
\end{proof}
\end{prop}

Theorem \ref{s4: thm 2.1f} follows from Proposition \ref{s4: prop 4.1f} and the following result.

\begin{lem}\label{s4: lem 4.2f}
Assume that $\Omega$ and $u$ satisfy the hypotheses of Theorem \ref{s4: thm 2.1f}. Then for any $h\le c(\rz)$, we have
\begin{equation}\label{s4: eq 4.5f}
d_n\ge ch^{\f{1}{2-\az}}.
\end{equation}
\end{lem}

$\mathit{Lemma\;\ref{s4: lem 4.2f}\;implies\;Theorem\;\ref{s4: thm 2.1f}}$.

From Lemma \ref{s4: lem 4.2f} and Proposition \ref{s4: prop 4.1f} we obtain
$$ch^{\f{1}{2}}\le d_i\le Ch^{\f{1}{2}},\;i=1,\dots,n-1,\quad\quad ch^{\f{1}{2-\az}}\le d_n\le Ch^{\f{1}{2-\az}}.$$
It follows that
\begin{equation}\label{s4: eq *7f}
\tilde{x}^*_h+cF_h B_1\subset A_h S_h\subset CF_h B_1,
\end{equation}
where we recall from Section 2 that
$$F_h x=(h^{\f{1}{2}}x', h^{\f{1}{2-\az}}x_n).$$
Since the domain of definition of $\tilde{G}_h$ contains a ball of radius $(\mu h/2)^{\f{1}{2}}$, we have
\begin{equation}\label{s4: eq *8f}
cF_h B_1\cap A_h\overline{\Omega}\subset A_h S_h\subset CF_h B_1.
\end{equation}
It follows that
\begin{equation}\label{s4: eq E_h'f}
c\mathcal{E}_h\cap A_h\overline{\Omega}\subset A_h S_h\subset C\mathcal{E}_h.
\end{equation}

Denote $A_h x=x-\nu_h x_n$. Using in \eqref{s4: eq *7f} that $S_{h/2}\subset S_h$ we find
$$|\nu_{h/2}-\nu_h|\le Ch^{\f{1}{2}-\f{1}{2-\az}},\quad\forall h\le c(\rz),$$
which gives
\begin{eqnarray}\label{s4: eq nu_hf}
|\nu_h|\le C(\rz)h^{\f{1}{2}-\f{1}{2-\az}},\quad\forall h\le c(\rz).
\end{eqnarray}
This easily implies that
\begin{equation}\label{s4: eq E_hf}
\mathcal{E}_{c_1(\rz)h}\subset A_h^{-1}\mathcal{E}_h\subset\mathcal{E}_{C_1(\rz)h}
\end{equation}
for some constants $c_1(\rz), C_1(\rz)>0$.

The conclusion of Theorem \ref{s4: thm 2.1f} follows from \eqref{s4: eq E_h'f} and \eqref{s4: eq E_hf}.

\bigskip

In order to prove Lemma \ref{s4: lem 4.2f}, we modify the definition of the quantity $b_u(h)$ in \cite{S2}. 

Fix $\az\in(0,1)$. Given a convex function $u$ we define
\begin{equation}\label{s4: b_u(h)f}
b_u(h)=h^{-\f{1}{2-\az}}\sup_{S_h}x_n.
\end{equation}
Whenever there is no possibility of confusion we drop the subindex $u$ and write for simplicity $b(h)$. 

$b(h)$ satisfies the following properties which are slightly different from those in \cite{S2}.
\begin{enumerate}
\item[1)] If $h_1\le h_2$, then
$$\l(\f{h_1}{h_2}\r)^{\f{1-\az}{2-\az}}\le\f{b(h_1)}{b(h_2)}\le\l(\f{h_2}{h_1}\r)^{\f{1}{2-\az}}.$$
\item[2)] If $A$ is a linear transformation which leaves the $x_n$-coordinate invariant and
$$\tilde{u}(Ax)=u(x),$$
then
$$b_{\tilde{u}}(h)=b_u(h).$$
\item[3)] If $A$ is a linear transformation which leaves the plane $\{x_n=0\}$ invariant, then
$$\f{b_{\tilde{u}}(h_1)}{b_{\tilde{u}}(h_2)}=\f{b_u(h_1)}{b_u(h_2)}.$$
\item[4)] If 
$$\tilde{u}(x)=\bz u(x)$$
with $\bz$ a positive constant, then
$$b_{\tilde{u}}(\bz h)=\bz^{-\f{1}{2-\az}}b_u(h)$$
and therefore
$$\f{b_{\tilde{u}}(\bz h_1)}{b_{\tilde{u}}(\bz h_2)}=\f{b_u(h_1)}{b_u(h_2)}.$$
\end{enumerate}

From part $(ii)$ of Proposition \ref{s4: prop 4.1f} we know that
$$cd_n\le d_h=x^*_h\cdot e_n\le Cd_n,$$
and it follows that
$$cd_n\le b_u(h)h^{\f{1}{2-\az}}=\sup_{S_h}x_n\le C d_n.$$
Thus Lemma \ref{s4: lem 4.2f} will follow if we show that $b_u(h)$ is bounded below. This will follow from property $1)$ above and the following lemma.

\begin{lem}\label{s4: lem 4.3f}
If $h\le c(\rz)$ and $b_u(h)\le c_0$, then
$$\f{b_u(th)}{b_u(h)}>2$$
for some $t\in[c_0,1]$.
\end{lem}

In order to prove Lemma \ref{s4: lem 4.3f}, we recall the function $\tilde{u}$, the section $\tilde{S}_h$ and the matrix $D_h$ in Proposition \ref{s4: prop 4.1f}. Define
$$v(x)=\f{1}{h}\tilde{u}(D_h x)=\f{1}{h}u(A_h^{-1}D_h x).$$
The section $S_1(v)=\{v<1\}=D_h^{-1}A_h S_h$ satisfies
$$x^*+cB_1\subset S_1(v)\subset CB_1$$
with $x^*$ the center of mass of $S_1(v)$.
The function $v$ satisfies in $S_1(v)$
$$\lz x_n^{-\az}\le\mathrm{det}\;D^2 v(x)=d_n^\az\mathrm{det}\;D^2 u(A_h^{-1}D_h x)\le\Lz x_n^{-\az}$$
and
$$v(0)=0,\quad 0\le v\le 1.$$

Moreover, let $0<t\le1$, $x^*_t(v)$ and $x^*_{th}$ be the centers of mass of $S_t(v)$ and $S_{th}(u)$ respectively, and $d_{th}=x^*_{th}\cdot e_n$. Then
$$(x^*_t(v)\cdot e_n)^{-\az}|S_t(v)|^2=\f{d_{th}^{-\az}|S_{th}(u)|^2}{h^n}.$$
Since Proposition \ref{s4: prop 4.1f} implies that $c(th)^n\le d_{th}^{-\az}|S_{th}(u)|^2\le C(th)^n$, we obtain
$$c t^n\le (x^*_t(v)\cdot e_n)^{-\az}|S_t(v)|^2\le C t^n.$$
From the convexity of $u$ we have
$$x^*_t(v)\cdot e_n=\f{d_{th}}{d_n}\ge c\f{d_{th}}{d_h}\ge c\cdot\f{\sup_{S_{th}(u)}x_n}{\sup_{S_h(u)}x_n}\ge ct.$$

Let $G_v:=\partial S_1(v)\cap\{v<1\}$. We claim that
$$G_v\subset\{x_n\le\sz\},\quad\quad\sz=C(\rz)h^{\f{1-\az}{n+1-\az}}.$$
Indeed, for $x\in G_v=D_h^{-1}\tilde{G}_h$,
$$d_n x_n\le\f{2}{\rz}|D'_h x'|^2\le C(\rz)h,$$
which gives
$$x_n\le C(\rz)h^{1-\f{n}{n+1-\az}}=\sz$$
by \eqref{s4: eq 4.2f}. Thus the claim follows.

We also have
$$v=1\quad\mathrm{on}\;\partial S_1(v)\setminus G_v.$$
On $G_v$,
$$\mu\sum_1^{n-1}a_i^2 x_i^2\le v(x)\le\mu^{-1}\sum_1^{n-1}a_i^2 x_i^2,$$
where
$$a_i=\f{d_i}{h^{\f{1}{2}}}\ge c_1,\quad i=1,\dots, n-1$$
by \eqref{s4: eq *4f}. 

In order to prove Lemma \ref{s4: lem 4.3f}, we only need to show that there exist constants $c(\rz), c_0$ small and $C$ sufficiently large such that if $h\le c(\rz)$ and $\max_{1\le i\le n-1}a_i\ge C$, then the rescaled function $v$ satisfies
\begin{equation}\label{s4: eq *9f}
b_v(t)\ge 2 b_v(1)
\end{equation}
for some $t\in[c_0,1]$.

\section{Proof of Theorem \ref{s4: thm 2.1f} (II)}\label{s4}

We consider the class of solutions $v$ that satisfy the properties above. After relabeling the constants $\mu$ and $a_i$, and by abuse of notation writing $u$ instead of $v$, we may assume we are in the following case.

Fix $\mu,\lz,\Lz$ and $\az\in(0,1)$. For an increasing sequence
$$a_1\le a_2\le\cdots\le a_{n-1}$$
with
$$a_1\ge\mu,$$
we consider the family of solutions
$$u\in\mathcal{D}^\mu_\sz(a_1,a_2,\dots, a_{n-1})$$
of convex functions $u:\Omega\to\mathbb{R}$ that satisfy
\begin{equation}\label{s4: eq 5.1f}
\lz x_n^{-\az}\le\mathrm{det}\;D^2 u\le\Lz x_n^{-\az},\quad\quad 0\le v\le 1\quad\mathrm{in}\;\Omega;
\end{equation}
\begin{equation}\label{s4: eq 5.2f}
0\in\partial\Omega,\quad\quad B_\mu(x_0)\subset\Omega\subset B^+_{1/\mu};
\end{equation}
\begin{equation}\label{s4: eq 5.3f}
\mu h^n\le (x^*_h\cdot e_n)^{-\az}|S_h|^2\le\mu^{-1}h^n,\quad\quad x^*_h\cdot e_n\ge\mu h
\end{equation}
with $x^*_h$ the center of mass of $S_h$.

Moreover, there exists a closed set $G\subset\partial\Omega$ such that
\begin{equation}\label{s4: eq 5.4f}
G\subset\partial\Omega\cap\{x_n\le\sz\},
\end{equation}
and $G$ is a graph in the $e_n$ direction with projection $\pi(G)$ along $e_n$,
\begin{equation}\label{s4: eq 5.5f}
\{\mu^{-1}\sum_1^{n-1}a_i^2 x_i^2\le 1\}\subset\pi(G)\subset\{\mu\sum_1^{n-1}a_i^2 x_i^2\le 1\}.
\end{equation}
The boundary values of $u=\varphi$ on $\partial\Omega$ satisfy
\begin{equation}\label{s7: eq 5.6}
\varphi=1\quad\mathrm{on}\;\partial\Omega\setminus G,
\end{equation}
and
\begin{equation}\label{s4: eq 5.7f}
\mu\sum_1^{n-1}a_i^2 x_i^2\le \varphi(x)\le\min\{1,\mu^{-1}\sum_1^{n-1}a_i^2 x_i^2\}\quad\mathrm{on}\;G.
\end{equation}

As explained in \cite{S2} (see Page 79 there), Property \eqref{s4: eq *9f} is a corollary of the following proposition.

\begin{prop}\label{s4: prop 5.1f}
For any $M>0$ there exists $C_*$ depending only on $M, n,\mu,\lz,\Lz$ and $\az$ such that if $u\in\mathcal{D}^\mu_\sz(a_1,a_2,\dots,a_{n-1})$ with
$$a_{n-1}\ge C_*,\quad\quad\sz\le C_*^{-1}$$
then
$$b(h)=(\sup_{S_h}x_n)h^{-\f{1}{2-\az}}\ge M$$
for some $h\in[C_*^{-1},1]$.
\end{prop}

We prove Proposition \ref{s4: prop 5.1f} by compactness as in \cite{S2}. We introduce the limiting solutions from the class $\mathcal{D}^\mu_{\sz}(a_1,\dots, a_{n-1})$ when $a_{k+1}\to\infty$ and $\sz\to 0$.

For an increasing sequence
$$a_1\le a_2\le\cdots\le a_{k}$$
with
$$a_1\ge\mu,$$
we denote by 
$$\mathcal{D}^\mu_0(a_1,\dots,a_k,\infty,\dots,\infty),\quad 0\le k\le n-2,$$
the class of functions $u$ that satisfy 
\begin{equation}\label{s4: eq g**5.1}
\lz x_n^{-\az}\le\mathrm{det}\;D^2 u\le\Lz x_n^{-\az},\quad\quad 0\le u\le 1\quad\mathrm{in}\;\Omega;
\end{equation}
\begin{equation}\label{s4: eq g*5.2}
0\in\partial\Omega,\quad\quad B_\mu(x_0)\subset\Omega\subset B^+_{1/\mu};
\end{equation}
\begin{equation}\label{s4: eq g*5.3}
\mu h^n\le (x^*_h\cdot e_n)^{-\az}|S_h|^2\le\mu^{-1}h^n,\quad\quad x^*_h\cdot e_n\ge\mu h,
\end{equation}
where $x^*_h$ is the center of mass of $S_h$. There exists a closed set $G$ such that
\begin{equation}\label{s4: eq g*5.8}
G\subset\partial\Omega\cap\{x_i=0,i>k\}.
\end{equation}
If we restrict to the space generated by the first $k$ coordinates, then
\begin{equation}\label{s4: eq g*5.9}
\{\mu^{-1}\sum_1^{k}a_i^2 x_i^2\le 1\}\subset G\subset\{\mu\sum_1^{k}a_i^2 x_i^2\le 1\}.
\end{equation}
The boundary values of $u=\varphi$ on $\partial\Omega$ satisfy
\begin{equation}\label{s4: eq g*5.10}
\varphi=1\quad\mathrm{on}\;\partial\Omega\setminus G,
\end{equation}
and
\begin{equation}\label{s4: eq g*5.11}
\mu\sum_1^{k}a_i^2 x_i^2\le \varphi(x)\le\min\{1,\mu^{-1}\sum_1^{k}a_i^2 x_i^2\}\quad\mathrm{on}\;G.
\end{equation}

As in \cite{S2}, Proposition \ref{s4: prop 5.1f} will follow from the proposition below.

\begin{prop}\label{s4: prop 5.2'}
For any $M>0$ and $0\le k\le n-2$ there exists $c_k$ depending only on $M, k, n,\mu,\lz,\Lz$ and $\az$ such that if 
$$u\in\mathcal{D}^\mu_0(a_1,\dots, a_{k}, \infty, \dots,\infty),$$ 
then
$$b(h)=(\sup_{S_h}x_n)h^{-\f{1}{2-\az}}\ge M$$
for some $h\in[c_k,1]$.
\end{prop}

To prove the above proposition, we use the notation introduced in \cite{S2}. Denote
$$x=(y, z, x_n),\quad y=(x_1,\dots, x_k)\in\mathbb{R}^k,\quad z=(x_{k+1},\dots, x_{n-1})\in\mathbb{R}^{n-1-k}.$$
A sliding along the $y$ direction is defined as follows:
$$T x:=x+\nu_1 z_1+\nu_2 z_2+\cdots+\nu_{n-k-1}z_{n-k-1}+\nu_{n-k}x_n$$
with 
$$\nu_1, \nu_2, \dots, \nu_{n-k}\in span\{e_1, \dots, e_k\}.$$

\begin{lem}\label{s3: lem 5.4}
Assume that
$$u\ge p(|z|-q x_n)$$
for some $p, q>0, q\le q_0$ and assume that for each section $S_h$ of $u$, $h\in (0,1)$, there exists $T_h$, a sliding along the $y$ direction, such that
$$T_h S_h\subset C_0 F_h B^+_1$$
for some constant $C_0$. Then
$$u\notin\mathcal{D}^\mu_0(1, \dots, 1, \infty, \dots, \infty).$$

\begin{proof}
Assume by contradiction that $u\in\mathcal{D}^\mu_0(1,\dots,1,\infty,\dots,\infty).$ We will show that
\begin{equation}\label{s3: eq 5.13}
u\ge p'(|z|-q' x_n),\quad\quad q'=q-\eta,
\end{equation}
where $\eta>0$ depends only on $q_0,C_0,\Lz,\mu, n, \az$ and $0<p'\ll p$.

Apply this result a finite number of times we obtain
$$u\ge\ez(|z|+x_n)$$
for some $\ez>0$ small. Thus we obtain $S_h\subset\{x_n\le\ez^{-1}h\}$ and it follows that
$$T_h S_h\subset\{x_n\le\ez^{-1}h\}.$$
This together with the hypothesis of the lemma and \eqref{s4: eq g*5.3} in the definition of the class $\mathcal{D}^\mu_0$ implies that 
$$\mu h^{n}\le (x^*_h\cdot e_n)^{-\az}|S_h|^2=(x^*_h\cdot e_n)^{-\az}|T_h S_h|^2\le Ch^{n+1-\az},$$
where $C$ is a constant depending only on $\ez, C_0, n, \mu$ and $\az$. This is a contradiction as $h\to 0$.

It remains to prove \eqref{s3: eq 5.13}. Since $u\in\mathcal{D}^\mu_0(1,\dots,1,\infty,\dots,\infty)$, there is a closed set 
$$G_h\subset\partial S_h\cap\{(z, x_n)=0\}$$ 
such that when we restrict to the subspace $\{(z, x_n)=0\}$, 
$$\{\mu^{-1}|y|^2\le h\}\subset G_h\subset\{\mu|y|^2\le h\},$$
and the boundary values $\varphi_h$ of $u$ on $\partial S_h$ satisfy
$$\varphi_h=h\quad\mathrm{on}\;\partial S_h\setminus G_h;$$
$$\mu|y|^2\le\varphi_h\le\min\{h,\mu^{-1}|y|^2\}\quad\mathrm{on}\;G_h.$$

Define
$$w(x)=\f{1}{h}u(T_h^{-1}F_h x).$$
Then
$$S_1(w)=F_h^{-1}T_h S_h\subset C_0 B^+_1,$$
and
$$\lz x_n^{-\az}\le\mathrm{det}\;D^2 w\le\Lz x_n^{-\az}\quad\mathrm{in}\;S_1(w).$$
Also,
\begin{equation}\label{s3: eq 5.14}
w(x)\ge\f{1}{h}p(|h^{\f{1}{2}}z|-q h^{\f{1}{2-\az}}x_n)=\f{p}{h^{\f{1}{2}}}(|z|-q h^{\f{\az}{2(2-\az)}}x_n).
\end{equation}

Moreover, the boundary values $\varphi_w$ of $w$ on $\partial S_1(w)$ satisfy
$$\varphi_w=1\quad\mathrm{on}\;\partial S_1(w)\setminus G_w;$$
$$\mu|y|^2\le\varphi_w\le\min\{1,\mu^{-1}|y|^2\}\quad\mathrm{on}\;G_w=F_h^{-1}G_h.$$
Define
$$v:=\dz\l(|x'|^2+\f{x_n^{2-\az}}{(2-\az)(1-\az)}\r)+\f{\Lz}{\dz^{n-1}}\l(z_1-q h^{\f{\az}{2(2-\az)}}x_n\r)^2+N\l(z_1-q h^{\f{\az}{2(2-\az)}}x_n\r)+\dz h^{\f{\az}{2(2-\az)}}x_n,$$
where $\dz$ is small depending only on $\mu, C_0, \az$ and $N$ is large such that
$$\f{\Lz}{\dz^{n-1}}t^2+Nt$$
is increasing in the interval $|t|\le(1+q_0)C_0$.

By straightforward computation and similar arguments to the proof of \cite[Lemma 5.4]{S2}, we find that $v$ is a lower barrier for $w$ in $S_1(w)$, which implies
$$w\ge N\l(z_1-q h^{\f{\az}{2(2-\az)}}x_n\r)+\dz h^{\f{\az}{2(2-\az)}}x_n\quad\mathrm{in}\;S_1(w).$$
Since this inequality holds for all directions in the $z$-plane, we obtain
$$w\ge N\l[|z|-\l(q-\f{\dz}{N}\r)h^{\f{\az}{2(2-\az)}}x_n\r].$$
Back to $u$ we have
\begin{eqnarray*}
u(x)=hw(F_h^{-1}T_h x)\ge h^{\f{1}{2}}N\l[|z|-\l(q-\f{\dz}{N}\r)x_n\r]\quad\mathrm{in}\;S_h.
\end{eqnarray*}
From the convexity of $u$ and $u(0)=0$, we know that this inequality holds in $\Omega$ and therefore \eqref{s3: eq 5.13} is proved.
\end{proof}
\end{lem}

Now we give the proof of Proposition \ref{s4: prop 5.2'}.

$\mathbf{k=0:}$ Assume Proposition \ref{s4: prop 5.2'} is not true, then by compactness, there exist $M>0$ and $u\in\mathcal{D}^\mu_0(\infty,\dots,\infty)$ such that $b(h)\le M$ for any $0<h\le 1$. Let
$$v:=\dz\l(|x'|+\f{1}{2}|x'|^2\r)+\f{\Lz}{\dz^{n-1}(2-\az)(1-\az)}x_n^{2-\az}-N x_n,$$ 
where $\dz$ is small depending only on $\mu$ and $N$ is large such that
$$\f{\Lz}{\dz^{n-1}(2-\az)(1-\az)}x_n^{2-\az}-N x_n\le 0$$
in $B^+_{1/\mu}$. It is easily seen that
$$v\le u\quad\quad\mathrm{in}\;\Omega.$$
It follows that
$$u\ge\dz|x'|-N x_n$$
and then
$$S_h\subset\{|x'|\le\dz^{-1}(N x_n+h)\}.$$
Since $b(h)\le M$ implies that $x_n\le M h^{\f{1}{2-\az}}\le M h^{\f{1}{2}}$, we obtain
$$S_h\subset\{|x'|\le Ch^{\f{1}{2}}, x_n\le Mh^{\f{1}{2-\az}}\},$$
where $C$ is a constant depending only on $M,\mu,\Lz$ and $\az$. This contradicts Lemma \ref{s3: lem 5.4} and therefore Proposition \ref{s4: prop 5.2'} is true for $k=0$. 

$\mathbf{Assume\;Proposition\;\ref{s4: prop 5.2'}\;holds\;for\;0, 1,\dots, k-1,\;1\le k\le n-2,\;and\; now\;we\;prove\;it\;for\;k.}$ 

By the induction hypothesis, it suffices to consider the case $a_k\le C_k$, where $C_k$ is a constant depending only on $M, k, n, \mu, \lz, \Lz$ and $\az$. Assume in contradiction that no $c_k$ exists, then we can find a limiting solution $u$ such that
\begin{equation}\label{s3: eq 5.15}
u\in\mathcal{D}^{\tilde{\mu}}_0(1,\dots,1, \infty, \dots, \infty)
\end{equation}
with 
\begin{equation}\label{s3: eq 5.16}
b(h)\le M,\quad\forall h>0,
\end{equation}
where $\tilde{\mu}$ depends only on $\mu$ and $C_k$.

Denote as before 
$$x=(y, z, x_n),\quad y=(x_1,\dots, x_k)\in\mathbb{R}^k,\quad z=(x_{k+1},\dots, x_{n-1})\in\mathbb{R}^{n-1-k}.$$
Similar to the case $k=0$, the function 
$$v:=\dz\l(|z|+\f{1}{2}|x'|^2\r)+\f{\Lz}{\dz^{n-1}(2-\az)(1-\az)}x_n^{2-\az}-N x_n$$ 
is a lower barrier for $u$, where $\dz$ is small depending only on $\tilde{\mu}$ and $N$ is large. Therefore,
\begin{equation}\label{s3: eq 5.17}
u\ge\dz|z|-N x_n.
\end{equation}
This together with \eqref{s3: eq 5.16} implies that
\begin{equation}\label{s3: eq 5.18}
S_h\subset\{|z|\le\dz^{-1}(N x_n+h)\}\cap\{x_n\le Mh^{\f{1}{2-\az}}\}.
\end{equation}

From John's lemma, there is an ellipsoid $E_h$ such that
\begin{equation}\label{s3: eq 5.19}
E_h\subset S_h-x^*_h\subset C(n) E_h
\end{equation}
with $x^*_h$ the center of mass of $S_h$. By a fact in linear algebra (see the arguments in \cite[Page 83]{S2}), there is $T_h$, a sliding along the $y$ direction, such that
\begin{equation}\label{s3: eq 5.20}
T_h E_h=|E_h|^{\f{1}{n}}A B_1,
\end{equation}
where, after rotating coordinates in the $(y,0,0)$ and $(0,z,0)$ subspaces, the matrix $A$ satisfies
$$A(y,z,x_n)=(A_1y, A_2(z,x_n)),$$
$$
A_1=\l(
\begin{array}{cccc}
\bz_1&0&\dots&0\\
0&\bz_2&\dots&0\\
\vdots&\vdots&\ddots&\vdots\\
0&0&\dots&\bz_k
\end{array}\r)\quad\mathrm{and}\quad A_2=
\l(\begin{array}{ccccc}
\gz_{k+1}&0&\dots&0&\tz_{k+1}\\
0&\gz_{k+2}&\dots&0&\tz_{k+2}\\
\vdots&\vdots&\ddots&\vdots&\vdots\\
0&0&\dots&\gz_{n-1}&\tz_{n-1}\\
0&0&\dots&0&\tz_n
\end{array}\r)
$$
with 
$$
0<\bz_1\le\dots\le\bz_k,\quad\quad\gz_j>0,\quad \tz_n>0,\quad\quad\l(\prod_1^{k}\bz_i\r)\l(\prod_{k+1}^{n-1}\gz_j\r)\tz_n=1.$$

Let 
$$\tilde{u}(x)=u(T_h^{-1}x),\quad\quad\tilde{S}_h=T_h S_h,$$
then \eqref{s3: eq 5.19} implies that
\begin{equation}\label{s3: eq 5.21}
\tilde{x}^*_h+|E_h|^{\f{1}{n}}A B_1\subset\tilde{S}_h\subset C(n)|E_h|^{\f{1}{n}}A B_1,
\end{equation}
where $\tilde{x}^*_h$ is the center of mass of $\tilde{S}_h$.

Since $u\in\mathcal{D}^{\tilde{\mu}}_0(1,\dots,1,\infty,\dots,\infty)$, there exists $\tilde{G}_h=G_h$,
$$\tilde{G}_h\subset\{(z, x_n)=0\}\cap\partial\tilde{S}_h$$
such that on the subspace $\{(z, x_n)=0\}$,
$$\{\tilde{\mu}^{-1}|y|^2\le h\}\subset\tilde{G}_h\subset\{\tilde{\mu}|y|^2\le h\},$$
and the boundary values $\tilde{\varphi}_h$ of $\tilde{u}$ on $\partial\tilde{S}_h$ satisfy
$$\tilde{\varphi}_h=h\quad\mathrm{on}\;\partial\tilde{S}_h\setminus\tilde{G}_h;$$
$$\tilde{\mu}|y|^2\le\tilde{\varphi}_h\le\min\{h,\tilde{\mu}^{-1}|y|^2\}\quad\mathrm{on}\;\tilde{G}_h.$$
For any $h>0$, denote $d_h:=x^*_h\cdot e_n$, then
\begin{equation}\label{s3: eq 6}
\tilde{\mu}h^n\le(\tilde{x}^*_h\cdot e_n)^{-\az}|\tilde{S}_h|^2=d_h^{-\az}|S_h|^2\le\tilde{\mu}^{-1}h^n,\quad\quad\tilde{x}^*_h\cdot e_n=d_h\ge\tilde{\mu}h.
\end{equation}

We will show that
\begin{equation}\label{s3: eq 7}
|E_h|^{\f{1}{n}}A B_1\subset C diag\l(h^{\f{1}{2}},\dots, h^{\f{1}{2}}, h^{\f{1}{2}}d_h^{\f{\az}{2}}\r) B_1,
\end{equation}
where $C$ is a constant depending only on $\mu, M, k, \lz, \Lz, n$ and $\az$.

This together with \eqref{s3: eq 5.21}, \eqref{s3: eq 5.16} gives
\begin{equation}\label{s3: eq 7'}
T_h S_h\subset C diag\l(h^{\f{1}{2}},\dots, h^{\f{1}{2}}, h^{\f{1}{2-\az}}\r) B_1.
\end{equation}

Now we prove \eqref{s3: eq 7}. Let 
\begin{eqnarray*}
\bar{A}&=&|E_h|^{\f{1}{n}}diag\l(h^{-\f{1}{2}},\dots, h^{-\f{1}{2}}, h^{-\f{1}{2}}d_h^{-\f{\az}{2}}\r) A\nonumber\\
&=&|E_h|^{\f{1}{n}}h^{-\f{1}{2}}\l(
\begin{array}{cc}
A_1&\\
&\l(
\begin{array}{cccc}
1&&&\\
&\ddots&&\\
&&1&\\
&&&d_h^{-\f{\az}{2}}
\end{array}\r)A_2
\end{array}\r)=\l(
\begin{array}{cc}
\bar{A}_1&\\
&\bar{A}_2
\end{array}\r).\nonumber\\
\end{eqnarray*}

Since $\tilde{G}_h\subset\partial\tilde{S}_h\cap\{(z, x_n)=0\}$ contains a ball in $\mathbb{R}^k$ of radius $(\tilde{\mu}h)^{1/2}$, then from the second inclusion in \eqref{s3: eq 5.21} we obtain
\begin{equation}\label{s3: eq 9}
\bar{\bz}_i:=h^{-\f{1}{2}}|E_h|^{\f{1}{n}}\bz_i\ge c,\quad i=1,\dots,k,
\end{equation}
where $c$ is a constant depending only on $n$ and $\tilde{\mu}$.

From \eqref{s3: eq 5.16} we know that for any $x=(y, z, x_n)\in\tilde{S}_h$ we have
\begin{eqnarray*}
x_n\le C(n)d_h&\le&C(n)(Mh^{\f{1}{2-\az}})^{\f{2-\az}{2}}d_h^{\f{\az}{2}}\le C(n,M,\az)h^{\f{1}{2}}d_h^{\f{\az}{2}},
\end{eqnarray*}
combining this, \eqref{s3: eq 6} and \eqref{s3: eq 5.18} we obtain that
$$\tilde{S}_h\subset\{|(z,x_n)|\le C h^{\f{1}{2}}d_h^{\f{\az}{2}}\}
.$$
This together with the first inclusion in \eqref{s3: eq 5.21} implies that $\|\bar{A}_2\|\le C$ and if follows that
\begin{equation}\label{s3: eq 10}
\bar{\gz}_j:=h^{-\f{1}{2}}|E_h|^{\f{1}{n}}\gz_j,\le C,\quad\quad h^{-\f{1}{2}}|E_h|^{\f{1}{n}}|\tz_n\nu|\le C,\quad\quad\bar{\tz}_n:=h^{-\f{1}{2}}|E_h|^{\f{1}{n}}d_h^{-\f{\az}{2}}\tz_n\le C,
\end{equation}
where $C$ is a constant depending only on $n,\tilde{\mu},\Lz,\az$ and $M$. 

Also, we have by \eqref{s3: eq 5.21} 
\begin{equation}\label{s3: eq 11}
|E_h|^{\f{1}{n}}\tz_n\le\tilde{x}^*_h\cdot e_n=d_h\le C(n)|E_h|^{\f{1}{n}}\tz_n.
\end{equation}

We define
$$w(x):=\f{1}{h}\tilde{u}(|E_h|^{\f{1}{n}}A x),$$
then from \eqref{s3: eq 5.21} we know that
$$B_1(x_0)\subset S_{1}(w)=|E_h|^{-\f{1}{n}}A^{-1}\tilde{S}_h\subset C(n)B_1$$
for some $x_0$, and from \eqref{s3: eq 6} and \eqref{s3: eq 11} we find that
$$\bar{\lz}x_n^{-\az}\le\mathrm{det}\;D^2 w\le\bar{\Lz}x_n^{-\az}$$
with $\bar{\lz},\bar{\Lz}$ depending only on $\lz, \Lz, n,\az,\tilde{\mu}$. 

Moreover, for $t>0$ let $x^*_t(w)$ be the center of mass of the section $S_t(w)$, then
$$S_t(w)=|E_h|^{-\f{1}{n}}A^{-1}T_h S_{th}(u),$$
and we have by \eqref{s3: eq 11}
$$\f{d_{th}}{d_h}\le x^*_t(w)\cdot e_n=|E_h|^{-\f{1}{n}}\tz_n^{-1}d_{th}\le C(n)\f{d_{th}}{d_h}.$$
Then \eqref{s3: eq 6} implies that
$$ct^n\le(x^*_t(w)\cdot e_n)^{-\az}|S_t(w)|^2\le C t^n$$
for some constants $c, C$ depending only on $n,\az$ and $\tilde{\mu}$. 

Let $G_w=\partial S_{1}(w)\cap\l\{w<1\r\}=|E_h|^{-\f{1}{n}}A^{-1}\tilde{G}_h$, then the boundary values $\varphi_w$ of $w$ satisfy
$$\varphi_w=1\quad\mathrm{on}\;\partial S_{1}(w)\setminus G_w,$$
and from the definition of $\bar{\bz}_i$ we find that
$$\tilde{\mu}\sum_1^k\bar{\bz}_i^2 y_i^2\le\varphi_w\le\tilde{\mu}^{-1}\sum_1^k\bar{\bz}_i^2 y_i^2.$$
This implies that
$$w\in\mathcal{D}^{\bar{\mu}}_0(\bar{\bz}_1,\dots,\bar{\bz}_k,\infty,\dots,\infty)$$
for some $\bar{\mu}$ depending only on $\mu, M, k, \lz, \Lz, n$ and $\az$.

Note that \eqref{s3: eq 6} implies that
\begin{eqnarray}\label{s3: eq 12}
c\le\l(\prod_1^k\bar{\bz}_i\r)\l(\prod_{k+1}^{n-1}\bar{\gz}_j\r)\bar{\tz}_n
=h^{-\f{n}{2}}|E_h|d_h^{-\f{\az}{2}}\le C
\end{eqnarray}
with $c, C$ depending only on $n$ and $\tilde{\mu}$. 

We claim
\begin{eqnarray}\label{s3: eq 13}
\bar{\tz}_n\ge c_*
\end{eqnarray}
for some small $c_*$ to be chosen.

Indeed, if we $c_*$ is small, then \eqref{s3: eq 10} and \eqref{s3: eq 12} imply that
$$\bar{\bz}_k\ge C_k(\bar{\mu},\bar{M},\bar{\lz},\bar{\Lz},n, \az)$$
with $\bar{M}:=2\bar{\mu}^{-1}$. Then by the induction hypothesis,
$$b_w(\bar{h})\ge\bar{M}\ge 2b_w(1)$$
for some $\bar{h}>C_k^{-1}$. It follows that
$$\f{b_u(h\bar{h})}{b_u(h)}=\f{b_w(\bar{h})}{b_w(1)}\ge 2,$$
which implies $b_u(h\bar{h})\ge 2 b_u(h)$ for any $h>0$. This contradicts \eqref{s3: eq 5.16} and therefore the claim holds.

Similarly, we obtain that
\begin{eqnarray}\label{s3: eq 14}
\bar{\gz}_j\ge \tilde{c}_*
\end{eqnarray}
for some small $\tilde{c}_*$.

We obtain from \eqref{s3: eq 12}, \eqref{s3: eq 13}, \eqref{s3: eq 14} that
\begin{eqnarray}\label{s3: eq 15}
\bar{\bz}_i\le C,\quad i=1,\dots, k,
\end{eqnarray}
where $C$ is a constant depending only on $\mu, M, k, \lz, \Lz, n$ and $\az$.
This implies that $\|\bar{A}_1\|\le C$ and therefore $\|\bar{A}\|\le C$. 

Thus, the estimate \eqref{s3: eq 7} holds. Then the proof is finished because \eqref{s3: eq 5.17}, \eqref{s3: eq 7'} and \eqref{s3: eq 5.15} contradict Lemma \ref{s3: lem 5.4}.

\section{Proof of Theorem \ref{s1: thm 2}}\label{s5}

In this section, we always denote by $c, C, c_i, C_i, i\in\mathbb{N}$ constants depending only on $n, c_0$ and $\az$ ($c_0$ is the constant in \eqref{s5: eq 7.3}). Their values may change from line to line whenever there is no possibility of confusion.

\begin{lem}\label{s5: lem 1}
Assume the hypotheses in Theorem \ref{s1: thm 2} hold, then for $i=1,\dots, n-1$ we have $u_i\in C(\overline{\mathbb{R}^n_+})$.

\begin{proof}
We first claim that for some constant $c_1$ small, we have
\begin{equation}\label{s5: lem 11}
|\nabla u|\le c_1^{-1},\quad\mathrm{in}\;B^+_{c_1}.
\end{equation}

Indeed, we note that 
\eqref{s5: eq 7.3} implies that
\begin{equation*}
B^+_k\subset S_1(u)\subset B^+_{k^{-1}}
\end{equation*}
for some constant $k$ depending only on $c_0$ and $\az$. We can use the convexity of $u$ and obtain an upper bound for $u_n$ and all $|u_i|, 1\le i\le n-1$, in $B^+_{k/4}$. On the other hand, for any $x_0\in B^+_{c_1}$, the function 
$$w_{(x'_0,0)}(x):=\f{1}{2}|x'_0|^2+x'_0\cdot(x'-x'_0)+\dz|x'-x'_0|^2+\f{\dz^{1-n}}{(2-\az)(1-\az)}(x_n^{2-\az}-k^{-1}x_n)$$
is a lower barrier $u$ in $S_1(u)$, where $\dz$ is small depending only on $n,c_0$ and $\az$. This together with the convexity of $u$ gives a lower bound for $u_n(x_0)$.

Next we prove that for any $1\le i\le n-1$, $u_i$ is continuous at any point $x_0\in\{|x'|\le c_1/2, x_n=0\}$.

Indeed, fix $1\le i\le n-1$ and $x_0\in\{|x'|\le c_1/2, x_n=0\}$, define
$$u_{x_0}(x)=u(x_0+x)-u(x_0)-\nabla u(x_0)\cdot x.$$
We only need to prove that $\partial_i u_{x_0}$ is continuous at $0$. Assume there is a sequence $x^{(m)}\to 0, m\to\infty$ with
$$\partial_i u_{x_0}(x^{(m)})\ge\ez$$
for some $\ez>0$. We have
$$u_{x_0}\ge u_{x_0}(x^{(m)})+\nabla u_{x_0}(x^{(m)})\cdot(x-x^{(m)}).$$
Note that $|\nabla u_{x_0}(x^{(m)})|$ is bounded by \eqref{s5: lem 11}. Let $m\to\infty$, we obtain 
$$u_{x_0}\ge a\cdot x$$
for some $a=(a', a_n)\in\mathbb{R}^n$ with $a_i\ge\ez$. From the value of $u_{x_0}$ on the boundary $\{x_n=0\}$ we find that $a'=0$. This is a contradiction.

For any $\lz>0$, we define
$$u_\lz(y):=\f{1}{\lz}u(F_{\lz}y),$$
then $u_\lz$ satisfies \eqref{s5: eq 7.3} and \eqref{s5: eq 7.4}. The results above show that for any $1\le i\le n-1$, $\partial_i u_{\lz}$ is continuous on $\{|x'|\le c_1/2, x_n=0\}$. Therefore, $u_i$ is continuous on $F_{\lz}\{|x'|\le c_1/2, x_n=0\}$. Let $\lz\to\infty$ and we conclude that $u_i$ is continuous on $\{x_n=0\}$.
\end{proof}
\end{lem}

$\mathbf{Proof\;of\;Theorem\;\ref{s1: thm 2}:}$ As before we have
\begin{equation}\label{s5: thm21}
B^+_k\subset S_1(u)\subset B^+_{k^{-1}}
\end{equation}
and
\begin{equation}\label{s5: thm22}
|u_i|\le C\quad\mathrm{in}\;B^+_{k/4},\quad i=1,\dots,n-1,
\end{equation}
where $k$ is a constant depending only on $c_0$ and $\az$.

Let 
$$L\varphi:=\mathrm{tr}[(D^2 u)^{-1}D^2\varphi]$$
be the linearized Monge-Amp$\grave{e}$re operator for $u$. Then for $i=1,\dots, n-1$ we have
\begin{eqnarray*}
Lu_i&=&0,\quad\quad u_i=x_i\quad\mathrm{on}\;\{x_n=0\},\\
Lu&=&n,
\end{eqnarray*}
and if we define $P(x)=\dz|x'|^2+\dz^{1-n}\f{x_n^{2-\az}}{(2-\az)(1-\az)}$ with $\dz>0$ a small constant to be chosen, then
\begin{eqnarray*}
LP=\mathrm{tr}[(D^2 u)^{-1}D^2 P]\ge n[\mathrm{det}(D^2 u)^{-1}\mathrm{det}\;D^2 P]^{\f{1}{n}}>n.
\end{eqnarray*}

Let $\gz_1, \gz_2$ be large constants to be chosen and define
$$v^{\mp}(x):=x_i\pm\gz_1\l[\dz|x'|^2+\dz^{1-n}\l(\f{x_n^{2-\az}}{(2-\az)(1-\az)}-\gz_2 x_n\r)-u(x)\r].$$
We have
$$Lv^-=\gz_1[LP-Lu]>0.$$
On $\partial B^+_{k/4}\cap\{x_n=0\}$, we choose $\dz\le 1/2$ and obtain
\begin{equation*}
v^-=x_i+\gz_1\l[\dz|x'|^2-\f{1}{2}|x'|^2\r]\le x_i.
\end{equation*}
We choose $\gz_2$ large such that
\begin{equation*}
\f{x_n^{2-\az}}{(2-\az)(1-\az)}-\gz_2 x_n\le 0\quad\mathrm{in}\;B^+_{k/4}.
\end{equation*}
Then on $\partial B^+_{k/4}\cap\{x_n>0\}$, we use \eqref{s5: eq 7.3} and obtain
\begin{eqnarray*}
v^-\le x_i+\gz_1\l[\dz|x'|^2-c_0|x'|^2-c_0 x_n^{2-\az}\r]\le x_i-\f{\gz_1 c_0}{2}(|x'|^2+x_n^{2-\az})\le-C,
\end{eqnarray*}
where $C$ is the constant in \eqref{s5: thm22}, and we choose $\dz\le c_0/2$ and $\gz_1$ large. 

By Lemma \ref{s5: lem 1}, $u_i\in C(\overline{B^+_{k/4}})$ and therefore the maximum principle for linear elliptic equations implies that  
$$v^-\le u_i\quad\mathrm{in}\;B^+_{k/4}.$$
Similarly,
$$v^+\ge u_i\quad\mathrm{in}\;B^+_{k/4}.$$
Therefore,
\begin{equation}\label{s5: thm23}
|u_i-x_i|\le\gz_1[\dz^{1-n}\gz_2 x_n+u]\quad\mathrm{in}\;B^+_{k/4}.
\end{equation}
For any $\lz>0$, we define
$$u_\lz(y):=\f{1}{\lz}u(F_{\lz}y),$$
then $u_\lz$ satisfies \eqref{s5: eq 7.3} and \eqref{s5: eq 7.4}.

Apply \eqref{s5: thm23} with $u\rightsquigarrow u_\lz$ and we obtain
\begin{equation*}
|\partial_iu_\lz(y)-y_i|\le\gz_1[\dz^{1-n}\gz_2 y_n+u_\lz(y)]\quad\quad\mathrm{in}\;B^+_{k/4}.
\end{equation*}
Back to $u$ we have
\begin{eqnarray}\label{s5: thm24}
|u_i(x)-x_i|\le\gz_1\l[\dz^{1-n}\gz_2\lz^{\f{1}{2}-\f{1}{2-\az}}x_n+\lz^{-\f{1}{2}}u(x)\r]\quad\mathrm{in}\;F_\lz B^+_{k/4}.
\end{eqnarray}
Let $\lz\to\infty$, we obtain
\begin{equation}\label{s5: thm25}
u_i=x_i,\quad\forall x\in\mathbb{R}^n_+.
\end{equation}
For any $x=(x',x_n)\in\mathbb{R}^n_+$,
\begin{eqnarray}\label{s5: thm26}
u(x', x_n)=u(0, x_n)+\int_0^1\nabla_{x'}u(\tz x', x_n)\cdot x'd\tz=u(0, x_n)+\f{1}{2}|x'|^2.
\end{eqnarray}
Thus, 
$$\mathrm{det}\;D^2 u=u_{nn}(0, x_n)=x_n^{-\az},$$
it follows that
$$u(0, x_n)=\f{x_n^{2-\az}}{(2-\az)(1-\az)}+A x_n+B$$
for some constants $A, B\in\mathbb{R}$.

Since $u\in C(\overline{\mathbb{R}^n_+})$ and satisfies \eqref{s5: eq 7.3}, we obtain $A=B=0$. The conclusion of the theorem follows from this and \eqref{s5: thm26}.

\section{Proof of Theorem \ref{s1: thm 3}}\label{s6}

$\mathbf{Proof\;of\;Theorem\;\ref{s6: thm 2.4}:\;}$ By the localization theorem (Theorem \ref{s4: thm 2.1}), 
\begin{equation*}
cU_0(x)\le u(x)\le c^{-1}U_0(x)\quad\mathrm{in}\;\Omega\cap S_{c}(U_0),
\end{equation*}
where $c$ is a constant depending only on $n,\az$ and $\rz$. Let $\Omega_h=F_h^{-1}\Omega$, then $\Omega_h\cap S_1(U_0)$ can be denoted by $x_n=g_h(x')$, where  
\begin{equation}\label{s6: thm 2.45}
g_h(x')=h^{-\f{1}{2-\az}}g(h^{\f{1}{2}}x')\le Ch^{\f{1-\az}{2-\az}}|x'|^2
\end{equation}
for some constant $C=C(n,\rz)$. Define
$$u_h(x)=\f{1}{h}u(F_h x),\quad x\in\Omega_h.$$
Then we have
\begin{equation}\label{s6: thm 2.46}
cU_0(x)\le u_h(x)\le c^{-1}U_0(x)\quad\mathrm{in}\;\Omega_h\cap S_1(U_0).
\end{equation}
The assumptions of Theorem \ref{s6: thm 2.4} imply that
\begin{equation}\label{s6: thm 2.47}
(1-\ez_0)(x_n-g_h(x'))^{-\az}\le\mathrm{det}\;D^2 u_h(x)\le (1+\ez_0)(x_n-g_h(x'))^{-\az}\quad\mathrm{in}\;\Omega_h\cap S_1(U_0)
\end{equation}
and
\begin{equation}\label{s6: thm 2.48}
\l(\f{1}{2}-\ez_0\r)|x'|^2\le u_h(x)\le\l(\f{1}{2}+\ez_0\r)|x'|^2\quad\mathrm{on}\;\partial\Omega_h\cap S_1(U_0).
\end{equation}

Assume by contradiction that Theorem \ref{s6: thm 2.4} does not hold. Then there is a constant $\eta>0$ such that for any $m\in\mathbb{N}, m\ge 1$, there exist $\Omega^m, g^m, u^m$ that satisfy the hypotheses of Theorem \ref{s6: thm 2.4} with $\ez_0\rightsquigarrow 1/m$, and some $0<h_m\le\min\{1/m, c\}$ such that if we denote $\Omega^m_{h_m}=F_{h_m}^{-1}\Omega^m$ and
$$u^m_{h_m}(x):=\f{1}{h_m}u^m(F_{h_m}x),$$
then the part $\partial\Omega^m_{h_m}\cap S_1(U_0)$ is given by $x_n=g^m_{h_m}(x')$ for some convex function $g^m_{h_m}$ satisfying \eqref{s6: thm 2.45} with $h\rightsquigarrow h_m$, and the function $u^m_{h_m}$ satisfies \eqref{s6: thm 2.46}-\eqref{s6: thm 2.48} with $\ez_0\rightsquigarrow 1/m$, while the inclusion in Theorem \ref{s6: thm 2.4} does not hold for $\eta$ and $S_{h_m}(u^m)$. 

Let $m\to\infty$, we can extract a subsequence $u^m_{h_m}$ that converges uniformly on compact sets to a global solution $u_0$ defined in $\mathbb{R}^n_+$ such that
\begin{equation}\label{s6: thm 2.49}
cU_0(x)\le u_0(x)\le c^{-1}U_0(x)\quad\mathrm{in}\;\mathbb{R}^n_+
\end{equation}
and
\begin{equation}\label{s6: thm 2.410}
\mathrm{det}\;D^2 u_0(x)=x_n^{-\az}\quad\mathrm{in}\;\mathbb{R}^n_+,\quad\quad u_0(x',0)=\f{1}{2}|x'|^2.
\end{equation}
Theorem 2 implies that
$$u_0=U_0=\f{1}{2}|x'|^2+\f{x_n^{2-\az}}{(2-\az)(1-\az)}.$$
We reach a contradiction.

$\mathbf{Proof\;of\;Theorem\;\ref{s1: thm 3}:\;}$ Assume the hypotheses in Theorem \ref{s1: thm 3} hold, then we can assume that 
\begin{equation*}
\varphi=\f{1}{2}\langle Mx', x'\rangle+o(|x'|^2),
\end{equation*}
for some positive definite matrix $M\in\mathbb{R}^{(n-1)\times(n-1)}$.

It suffices to prove the theorem for the case $f(0)=1$ and $M=I_{n-1}$. Indeed, let $D'\in\mathbb{R}^{(n-1)\times(n-1)}$ be a positive definite matrix such that
$$D'MD'=I_{n-1}.$$
Let $\lz>0$ be a constant to be chosen. Define
$$D:=\l(
\begin{array}{cc}
D'&\\
&\lz
\end{array}\r).$$
For any $y\in\tilde{\Omega}=D^{-1}\Omega$, define 
$$\tilde{u}(y)=u(Dy)$$
and $\tilde{\varphi}:=\tilde{u}|_{\partial\tilde{\Omega}}$. Then we have
\begin{eqnarray*}
\mathrm{det}\;D^2\tilde{u}(y)=\tilde{f}(y)d_{\partial\tilde{\Omega}}^{-\az}(y),\quad\quad\tilde{f}(y):=(\mathrm{det}\;M)^{-1}\lz^2 f(Dy)\f{d_{\partial\Omega}^{-\az}(Dy)}{d_{\partial\tilde{\Omega}}^{-\az}(y)}.
\end{eqnarray*}
It is easy to see that 
\begin{equation*}
\lim_{y\to 0}\f{d_{\partial\Omega}(Dy)}{d_{\partial\tilde{\Omega}}(y)}=\lz.
\end{equation*}
Thus we can choose $\lz>0$ such that 
$$\lim_{y\to 0}\tilde{f}(y)=(\mathrm{det}\;M)^{-1}\lz^2 f(0)\lz^{-\az}=1.$$

Now we assume $f(0)=1$ and $M=I_{n-1}$, and we will prove that
\begin{equation}\label{s6: eq *}
u(x)=\f{1}{2}|x'|^2+\f{x_n^{2-\az}}{(2-\az)(1-\az)}+o(|x'|^2+x_n^{2-\az}).
\end{equation}
For any $\ez_1>0$ small, we can choose $R=R(\ez_1)>0$ such that $\partial\Omega\cap B_R$ is given by $x_n=g(x')$ for some convex function $g$, where
\begin{equation}\label{s6: eq 3.1''}
g\in C^2\l(\pi(\partial\Omega\cap B_R)\r), \quad g(0)=0,\quad\nabla g(0)=0,\quad D^2 g(0)\ge k_0 I_{n-1}>0,
\end{equation}
\begin{equation}\label{s6: eq 3.2}
\mathrm{det}\;D^2 u=f(x)d_{\partial\Omega}^{-\az},\quad 1-\ez_1\le f\le 1+\ez_1\quad\mathrm{in}\;\Omega\cap B_R,
\end{equation}
\begin{equation}\label{s6: eq 3.4}
\l(\f{1}{2}-\ez_1\r)|x'|^2\le u(x)=\varphi(x')\le\l(\f{1}{2}+\ez_1\r)|x'|^2\quad\mathrm{on}\;\partial\Omega\cap B_R,
\end{equation}
where $k_0$ depends only on the principal curvatures of $\partial\Omega$ at $0$.

It is obvious that 
\begin{equation*}
\lim_{x\to 0}\f{x_n-g(x')}{d_{\partial\Omega}(x)}=1.
\end{equation*}
Therefore we can choose $R=R(\ez_1)$ smaller such that
\begin{equation}\label{s6: eq 3.2'}
(1-4\ez_1)(x_n-g(x'))^{-\az}\le\mathrm{det}\;D^2 u\le (1+4\ez_1)(x_n-g(x'))^{-\az}\quad\mathrm{in}\;\Omega\cap B_R,
\end{equation}

For any $\eta>0$, let $\ez_0$ be the constant given by Theorem \ref{s6: thm 2.4} and $\ez_1:=\ez_0/4$. Using \eqref{s6: eq 3.1''}, \eqref{s6: eq 3.4} and \eqref{s6: eq 3.2'}, we can choose $\rz>0$ such that the hypotheses of Theorem \ref{s6: thm 2.4} hold. Then Theorem \ref{s6: thm 2.4} implies that
$$|u(x)-U_0(x)|\le C\eta U _0(x)\quad\mathrm{in}\;\Omega\cap S_{c}(U_0),$$
where $U_0$ is defined as in Theorem \ref{s6: thm 2.4}, $c$ is a constant depending only on $\eta, n, \az, \rz$ and $C=C(n,\az)$ depends only on $n, \az$. This proves \eqref{s6: eq *} and therefore the proof of Theorem \ref{s1: thm 3} is complete.

\section{Proof of Theorem \ref{s1: thm 4}}\label{s7}
In this section we always denote by $c, C, c_i, C_i (i=0,1,\dots)$ constants depending only on $n,\lz_0,\Lz_0,\az$, $\mathrm{diam}(\Omega)$, and $\varphi,\partial\Omega$ up to their second derivatives. For any $A,B\in\mathbb{R}$, we write $A\sim B$ if 
$$c\le\f{A}{B}\le C$$
for some constants $c, C$ depending only on $n,\lz_0,\Lz_0,\az$, $\mathrm{diam}(\Omega)$, and $\varphi,\partial\Omega$ up to their second derivatives. 

Suppose the assumptions of Theorem \ref{s1: thm 4} hold. First we can use barriers to obtain that
\begin{equation}\label{s7: eq *0}
\|u\|_{C(\overline{\Omega})}\le C.
\end{equation}
Now we restrict to a neighborhood of $0\in\partial\Omega$. As in Section 2, we can assume that for some fixed small $\rz>0$, the part $\partial\Omega\cap\{x_n\le\rz\}$ is given by $x_n=g(x')$ for some convex function $g$, where
\begin{equation}\label{s7: eq 3.1''}
g\in C^2\l(\pi(\partial\Omega\cap\{x_n<\rz\})\r), \quad g(0)=0, \quad\nabla g(0)=0.
\end{equation}
The function $u:\overline{\Omega}\to\mathbb{R}$ satisfies $u=\varphi(x')$ on $\partial\Omega\cap\{x_n\le\rz\}$, and
\begin{equation}\label{s7: eq 3.2}
\mathrm{det}\;D^2 u=f,\quad 0<\lz(x_n-g)^{-\az}\le f\le\Lz(x_n-g)^{-\az}\quad\mathrm{in}\;\Omega\cap\{x_n<\rz/2\},
\end{equation}
where $\az\in[1,2)$.

\bigskip

$\mathbf{Case\;1:}$ $\az\in(1,2)$.

Denote $\bz:=\f{n+\az-1}{n}>1$. We claim that for any $x\in\partial\Omega\cap\{x_n\le\rz/2\}$,
\begin{equation}\label{s7: eq 3.4}
-\dz^{-1}x_n-\dz^{-1}(x_n-g)^{2-\bz}\le\tilde{u}\le\dz^{-1}x_n-\dz(x_n-g)^{2-\bz},
\end{equation}
where $\dz>0$ is a small constant, and 
$$\tilde{u}:=u-u(0)-\nabla_{x'}\varphi(0)\cdot x'.$$
Indeed, let $C_0, C_1>0$ be two constants and define
$$v^-:=\varphi(0)+\nabla_{x'}\varphi(0)\cdot x'-\f{C_0}{(2-\bz)}(x_n-g)^{2-\bz}-C_1x_n.$$
Since $\Omega$ is uniformly convex, $\partial\Omega,\varphi\in C^2$, and $u$ is bounded below by \eqref{s7: eq *0}, by straightforward computation we obtain that $v^-$ is a lower barrier for $u$ in $\Omega\cap\{x_n\le\rz/2\}$ if $C_0, C_1$ are sufficiently large. Similarly, the function
$$v^+:=\varphi(0)+\nabla_{x'}\varphi(0)\cdot x'-\f{c_1}{(2-\bz)}(x_n-g)^{2-\bz}+Cx_n$$
is an upper barrier for $u$ in $\Omega\cap\{x_n\le\rz/2\}$ if $c_1$ is small and $C$ is sufficiently large. Hence the claim follows.

The estimate \eqref{s7: eq 3.4} implies that
$$|u-u(x_0)|\le C|x-x_0|^{2-\bz}\quad\forall\;x_0\in\partial\Omega,\;x\in\Omega.$$
This together with the convexity of $u$ implies that $u$ is H$\ddot{o}$lder continuous in $\Omega$ and
\begin{equation}\label{s7: eq h}
\|u\|_{C^{0,2-\bz}(\overline{\Omega})}\le C.
\end{equation}

Let $y_0\in\Omega$ and assume $S_{\bar{h}}(y_0)$ is the maximal section included in $\Omega$ which becomes tangent to $\partial\Omega$ at $0$ with $\bar{h}\le c$.
Then we obtain that 
$$\nabla_{x'}u(y_0)=\nabla_{x'}\varphi(0),\quad\quad S_{\bar{h}}(y_0)=\{x\in\Omega: \tilde{u}(x)<u_n(y_0)x_n\}.$$ 
Note that $u_n(y_0)$ is bounded above since $\varphi\in C^2$ and $\Omega$ is uniformly convex. Thus 
$$\bar{h}=-\tilde{u}(y_0)+u_n(y_0)y_0\cdot e_n$$
is bounded above.

Denote $M:=-u_n(y_0)$. We only need to consider two cases: $-C<M<2\dz^{-1}$ and $M\ge 2\dz^{-1}$, where $\dz$ is the constant in \eqref{s7: eq 3.4}.

If $-C<M<2\dz^{-1}$, then at the point $x=(0, c_0)$ with $c_0$ a small constant, we have by \eqref{s7: eq 3.4}
$$\tilde{u}+Mx_n\le 3\dz^{-1}x_n-\dz x_n^{2-\bz}\le-\f{\dz}{2}x_n^{2-\bz}=-\f{\dz}{2}c_0^{2-\bz}.$$
It follows that $\bar{h}$ is bounded below. Hence by \eqref{s7: eq h},
$$S_{\bar{h}}(y_0)\supset B_c(y_0).$$

It remains to consider the case $M\ge 2\dz^{-1}$. For some $c_1$ small, the second inequality in \eqref{s7: eq 3.4} implies that the point $x=(0, c_1M^{\f{1}{1-\bz}})\in S_{\bar{h}}(y_0)$ and 
$$\tilde{u}+Mx_n\le-\f{\dz}{2}x_n^{2-\bz}+Mx_n=c_1M^{\f{2-\bz}{1-\bz}}\l[1-\f{\dz}{2}c_1^{1-\bz}\r]\le-c_1M^{\f{2-\bz}{1-\bz}}.$$
Therefore,
\begin{equation}\label{s7: eq *1}
\bar{h}\ge c_1M^{\f{2-\bz}{1-\bz}}.
\end{equation}
On the other hand, the first inequality in \eqref{s7: eq 3.4} and the uniform convexity of $\Omega$ imply that
\begin{equation}\label{s7: eq *2}
S_{\bar{h}}(y_0)\subset\{x_n\le CM^{\f{1}{1-\bz}}, |x'|\le CM^{\f{1}{2(1-\bz)}}\}.
\end{equation}
Using this and the first inequality in \eqref{s7: eq 3.4} again, we obtain that for any $x\in S_{\bar{h}}(y_0)$
$$\tilde{u}+Mx_n\ge\f{M}{2}x_n-\dz^{-1}x_n^{2-\bz}\ge-\dz^{-1}x_n^{2-\bz}\ge-CM^{\f{2-\bz}{1-\bz}}$$
and therefore
\begin{equation}\label{s7: eq *3}
\bar{h}\le CM^{\f{2-\bz}{1-\bz}}.
\end{equation}
By \eqref{s7: eq 3.4}, \eqref{s7: eq *1} and \eqref{s7: eq *3}, we have
\begin{equation}\label{s7: eq *4}
S_{\bar{h}/2}(y_0)\subset\{x_n\ge cM^{\f{1}{1-\bz}}\}.
\end{equation}
Using the first inequality in \eqref{s7: eq 3.4}, \eqref{s7: eq *2} and \eqref{s7: eq *4},we obtain that
\begin{equation}\label{s7: eq *4'}
x_n-g\sim x_n\sim M^{\f{1}{1-\bz}}\quad\quad\forall\;x\in S_{\bar{h}/2}(y_0).
\end{equation}
The volume estimate for interior sections in \cite[Corollary 3.2.4]{G} and the definition of $\bz$ imply that
\begin{equation}\label{s7: eq *5}
|S_{\bar{h}/2}(y_0)|\sim M^{\f{n(2-\bz)+\az}{2(1-\bz)}}=M^{\f{n+1}{2(1-\bz)}}.
\end{equation}
Define $Ty:=(M^{\f{1}{2(1-\bz)}}y', M^{\f{1}{1-\bz}}y_n)$, and
$$v(y):=\f{1}{M^{\f{2-\bz}{1-\bz}}}[\tilde{u}(Ty)+M(Ty)\cdot e_n+\bar{h}/2]\quad\mathrm{in}\;T^{-1}S_{\bar{h}/2}(y_0).$$
We have
$$\mathrm{det}\;D^2 v(y)=M^{\f{\az}{1-\bz}}\mathrm{det}\;D^2 u(Ty)\sim 1\quad\mathrm{in}\;T^{-1}S_{\bar{h}/2}(y_0)$$
and
$$|T^{-1}S_{\bar{h}/2}(y_0)|\sim 1.$$
Since \eqref{s7: eq *2} implies that
$$T^{-1}S_{\bar{h}/2}(y_0)\subset B_C,$$
hence the Aleksandrov's maximum principle \cite[Theorem 1.4.2]{G} implies that
$$T^{-1}S_{\bar{h}/2}(y_0)\supset B_c(T^{-1}y_0).$$
Part $i)$ of Theorem \ref{s1: thm 4} is proved.

\bigskip

$\mathbf{Case\;2:}$ $\az=1$.

By straightforward computation, the functions
$$v^-:=\varphi(0)+\nabla_{x'}\varphi(0)\cdot x'-C_0(x_n-g)(-\log(x_n-g))^{\f{1}{n}}-C_1x_n$$
and
$$v^+:=\varphi(0)+\nabla_{x'}\varphi(0)\cdot x'-c_1(x_n-g)(-\log(x_n-g))^{\f{1}{n}}+Cx_n,$$
are barriers for $u$ in $\Omega\cap\{x_n\le\rz/2\}$ if $C_0,C_1,C$ are large constants and $c_1$ is small. Hence, we have in $\Omega\cap\{x_n\le\rz/2\}$
\begin{equation}\label{s7: eq 3.4*}
-\dz^{-1}x_n-\dz^{-1}(x_n-g)\l(-\log(x_n-g)\r)^{\f{1}{n}}\le\tilde{u}\le\dz^{-1}x_n-\dz(x_n-g)\l(-\log(x_n-g)\r)^{\f{1}{n}},
\end{equation}
where $\dz>0$ is a small constant, and $\tilde{u}$ is defined as in the case $\az\in(1,2)$.

Using \eqref{s7: eq 3.4*} and similar arguments to the previous case, part $ii)$ of the theorem is proved. The proof of Theorem \ref{s1: thm 4} is complete.

\end{document}